\documentclass[10pt]{article}
\usepackage{amsmath,amssymb,wasysym}
\usepackage{graphicx}
\newtheorem{theo}{Theorem}[section]
\newtheorem{lemm}{Lemma}[section]
\newtheorem{coro}{Corolary}[section]
\newtheorem{prop}{Proposition}[section]
\newtheorem{defi}{Definition}[section]
\newtheorem{remark}{Remark}[section]
\def\proof {{\noindent \bf{Proof:\hspace{4pt}}}}
\def\endproof{\hfill$\square$\vspace{6pt}}
\numberwithin{equation}{section}
\title{
{\bf\Large  Nonlinear Stability of Periodic Travelling Wave Solutions for the Regularized Benjamin-Ono and BBM Equations }}
\author{{\bf\large Angulo, Jaime}\footnote{Email: angulo@ime.usp.br}\hspace{2mm}
{\bf\large}\vspace{1mm}\\
{\it\small Department of Mathematics, IME-USP}\\
 {\it\small Rua do Mat\~ao 1010, Cidade Universit\'aria,}\\
{\it\small  CEP 05508-090, S\~ao Paulo, SP, Brazil}\vspace{3mm}\\
{\bf\large Scialom, M\'arcia}\footnote{Email: scialom@ime.unicamp.br}\hspace{0.7cm}
{\bf\large Banquet, Carlos}\vspace{1mm}\footnote{Email: cbanquet@ime.unicamp.br}\\
{\it\small  Department of Mathematics, IMECC-UNICAMP}\\ 
{\it\small Rua S\'ergio Buarque de Holanda 651,}\\
 {\it\small CEP 13081-970, Campinas, SP, Brazil}}
\date{}
\begin{document}
\maketitle


\begin{abstract}
This paper has various goals: first, we develop a local and global well-posedness theory for the regularized Benjamin-Ono equation in the periodic setting, second, we show that the Cauchy problem for this equation (in both periodic and non-periodic case) cannot be solved by an iteration scheme based on the Duhamel formula for negative Sobolev indices, third, a proof of the existence of a smooth curve of periodic travelling wave solutions, for the regularized Benjamin-Ono equation, with fixed minimal period $2L$, is given. It is also shown that these solutions are nonlinearly stable in the energy space $H^{1/2}_{per}$ by perturbations of the same wavelength. Finally, an extension of the theory developed for the regularized Benjamin-Ono equation is given and as an example it is proved that the cnoidal wave solutions associated to the Benjamin-Bona-Mahony equation are nonlinearly stable in $H^1_{per}$.
\end{abstract}


\section{Introduction}
In this work we are interested in the study of the regularized Benjamin-Ono equation (rBO equation henceforth)
\begin{equation} \label{ecuabasica}
u_t+u_x+uu_{x}+\mathcal{H}u_{xt}=0,
\end{equation}
where $u$ is a real-valued function and $\mathcal{H}$ denotes the Hilbert transform defined via the Fourier transform as 
\[\widehat{\mathcal{H}f}(k) = -i sgn (k)\widehat{f}(k),\]
where
\[sgn(k)=\left \{ 
\begin{aligned}
-1 &, \ \ \ k<0\\
1 &, \ \ \ k>0.\\
\end{aligned} \right.\]
The regularized Benjamin-Ono equation is a model for the time evolution of long-crested waves at the interface between two immiscible fluids. Some situations in which the equation is useful are the pycnocline in the deep ocean, and the two-layer system created by the inflow of fresh water from a river into the sea, see \cite{kalisch1}. This equation is formally equivalent to the Benjamin-Ono equation 
\begin{equation}\label{boequation}
v_t+v_x+vv_x-\mathcal{H}v_{xx}=0,
\end{equation}
which was first introduced by Benjamin \cite{benjamin2} and later by Ono \cite{ono1}  as a model equation for the same situation as the rBO. Being more specific, for suitably restricted initial conditions, the solutions $u$ of (\ref{ecuabasica}) and $v$ of (\ref{boequation}) are nearly identical at least for values of $t$ in $[0,T]$ where $T$ is quite large, see \cite{bonaAlbert1} for more details. See also \cite{bonaAlbert1} for a more detail discussion about the advantages and disadvantages of using equation (\ref{boequation}) or (\ref{ecuabasica}) for modelling the propagation of small-amplitude long waves.\\

As far as we know in the periodic setting there does not exist any result about the well-posedness problem for the equation (\ref{ecuabasica}). We prove that the periodic initial value problem 
\begin{equation} \label{probperiod}
\left \{ 
\begin{aligned}
u_t+u&_x+uu_x+\mathcal{H}u_{xt}=0, \ x\in \mathbb{R}, \ t\in \mathbb{R},\\
u(x,0)&=u_0(x),\\
\end{aligned} \right.
\end{equation}
with initial data in periodic Sobolev spaces $H^s_{per}([-L,L])$ is locally wellposed if $s>1/2$ and globally wellposed if $s\geq 3/2.$ These results on existence, uniqueness and continuous dependence on the initial data of solutions to (\ref{probperiod}) have no special difficulty, and we were able to established them almost in the same way as in Bona and Kalisch \cite{kalisch2}, where they were obtained on the real line. Local well-posedness can be proved using a contraction argument in the space $C([0,T];H^s_{per})$ and  global solution is obtained via \textit{a priori} estimates.\\ 

The rBO equation possesses three conservation laws: 
\begin{equation}\label{prisegley}
E(u):= \frac{1}{2}\int (u\mathcal{H}u_x -\frac{1}{3} u^3)dx, \ \ \ \ \ \  F(u):= \frac{1}{2}\int (u^2+u\mathcal{H}u_{x})dx
\end{equation}
and 
\begin{equation*}
G(u):=\int udx.
\end{equation*}
This suggests that the space $H^{\frac 1{2}}(\mathbb{R})$ (or $H^{\frac{1}{2}}_{per}$) is a good candidate for a global well-posedness theory of the Cauchy problem associated to (\ref{ecuabasica}). This problem is open in  $H^s(\mathbb{R})$ (or $H^{s}_{per}$) with $s\leq 1/2$ and one of the goals of this paper is to present some obstructions to its solution by iteration methods. More precisely, we prove that the flow map data-solution cannot be $C^2$ for $s<0,$ in both, periodic and non-periodic case. This kind of ill-posedness was studied by Bourgain in \cite{Bourgain1} and Tzvetkov in \cite{Tzvetkov1} for the KdV equation; Molinet, Saut and Tzvetkov in \cite{MoliSautTzve1} and \cite{MoliSautTzve2} did the same for the Benjamin-Ono equation and the Kadomsev-Petviashvili I (KPI) equation, respectively. As far as we know, there does not exist any result about ill-posedness for the rBO on the periodic and non-periodic cases.\\

The last part of the paper is  dedicated to another important qualitative aspect of nonlinear dispersive equations, the travelling wave solutions, which depending on the specific boundary conditions on the wave's shape can be either solitary or periodic waves. The existence, nonlinear stability and instability of solitary wave solutions have been discussed in the past two decades from several points of view. Many techniques have been created to find solutions and sufficient conditions have been obtained to insure the stability or instability of this kind of waves, see for example \cite{Albert1}, \cite{bona1},\cite{bonaAlbert2}, \cite{benjamin1}, \cite{benjamin2}, \cite{bona2}, \cite{grillakis2}, \cite{grillakis3}, \cite{ono1}, \cite{weinstein2}, \cite{weinstein3}.
In contrast to the study of solitary waves, the periodic travelling wave solutions has received less attention. In recent years some papers in this subject have appeared, see for instance \cite{angulo3}, \cite{angulo4}, \cite{anguloBonaScialom}, \cite{anguloNatali}, \cite{GallayHaragus1}, \cite{GallayHaragus2}, \cite{NataliPastor}.\\

In this paper we are interested in giving a stability theory of periodic travelling wave solutions for the nonlinear dispersive equation (\ref{ecuabasica}). The periodic travelling wave solutions to be considered here will be of the general form
\[u(x,t)=\phi(x-ct),\]
where $\phi:\mathbb{R}\rightarrow \mathbb{R} $ is a smooth $2L-$periodic function and $c\neq 1.$ So, by replacing these permanent wave form into (\ref{ecuabasica}), integrating and considering the constant of integration equal to zero, we obtain
\begin{equation} \label{ecuaordina}
 c\mathcal{H}\phi_c' +(c-1)\phi_{c}-\frac{1}{2} \phi_{c}^2= 0.
\end{equation}
In the framework of travelling waves of solitary type, it is known the existence of solutions for (\ref{ecuaordina}) in the form
\begin{equation}\label{solitaria}
\phi_c(x)=\frac{4(c-1)}{1+\left(\frac{c-1}{c}x\right)^2}\ ,
\end{equation}
where $c>1.$ The stability theory for this solitary wave was established by Albert, Bona and Henry in \cite{bona1}. Additionally, Kalisch in \cite{kalisch1} exhibited a periodic family of travelling wave solutions (depending of the speed) with period $2\pi,$ for the rBO and used it to test the rate of convergence of a numerical scheme, which was introduced in \cite{kalisch2} to prove that the equation (\ref{ecuabasica}) does not constitute an \textit{infinite-dimensional completely integrable system.}\\

For $L>\pi$ and $c>1+\frac{\pi}{L-\pi},$ we prove the existence of a smooth curve of even periodic travelling wave solutions for (\ref{ecuabasica}). The construction of this solutions is based on the Poisson Summation theorem (as in Angulo and Natali \cite{anguloNatali}). The family of solutions is given by
\begin{equation}\label{solperiodicas}
\phi_{c}(\xi)= \frac{2c\pi}{L}\left( \frac{\sinh (\eta)}{\cosh(\eta)-\cos\left(\frac{\pi \xi}{L}\right)}\right),
\end{equation}
with $\eta$ satisfying
\begin{equation}\label{valortanh}
\eta(c)=\tanh^{-1}\left(\frac{c\pi}{(c-1)L}\right).
\end{equation}
Concerning the nonlinear stability of this family of periodic solutions we extend the classical approach developed by Benjamin \cite{benjamin1}, Bona \cite{bona1} and Weinstein \cite{weinstein2} to the periodic case. More precisely, we use the conservations laws  (\ref{prisegley}) to prove  that the orbit $\Omega_{\phi_c}=\left\{\phi_c(\cdot+y): y\in \mathbb{R} \right\}$ generated by the wave $\phi_c$ is orbitally stable in $H^{\frac{1}{2}}_{per}([-L,L])$ with regard to the periodic flow of the rBO equation.\\

In order to get the spectral conditions required to prove the nonlinear stability, we use the recent theory developed by Angulo and Natali \cite{anguloNatali}. Although their theory was established for another kind of equations we still can apply it to obtain an specific spectral structure associated to the non-local operator
\begin{equation}\label{operadorMio}
\mathcal{L}=c\mathcal{H}\partial_x-1+c-\phi_{c}.
\end{equation}

In the last section of the paper the theory established for the rBO equation is extended for a general family of regularized equations. We study a class of equations of the form
\begin{equation}\label{formagen}
u_t+u_x+u^pu_x+Hu_t=0,
\end{equation}
where $p\geq 1$ is an integer and $H$ is a differential or pseudo-differential operator in the context of periodic functions. Note that a considerable range of equations arise in practice. For instance, if you consider $H=-\partial_{x}^2$ you get the generalized Benjamin-Bona-Mahony equation and if $H=\mathcal{H}\partial_x$ we obtain the generalized regularized Benjamin-Ono equation, so the stability issues  for the generalized class (\ref{formagen}) are not just of mathematical interest.\\

We give sufficient conditions to get the nonlinear stability of periodic wave solutions associated to equations of the type (\ref{formagen}), and as example we prove that the cnoidal wave solutions of the BBM equation, with fundamental period $L>2\pi,$ are orbitally stable in $H^1_{per}([0,L]),$ for speeds $c>1+\frac{4\pi^2}{L^2-4\pi^2}.$ \\

This kind of generalization in the context of solitary waves have been studied before, see for example \cite{bona1} and \cite{bonaChen1}. In the periodic setting H\u{a}r\u{a}gu\c{s} in \cite{haragus1} proved the spectral stability of periodic travelling wave solutions, for the generalized BBM, which are small perturbations of the constant solution $u=(c-1)^{1/p},$ in both $L^2(\mathbb{R})$ and $C_b(\mathbb{R}).$ If $1\leq p\leq 2,$ she proved spectral stability for $c>1$ and for $p\geq 3,$ there exists a critical speed $c_p$ such that the periodic waves are spectrally stable for $c\in(c_p,\frac{p}{p-3})$, and unstable for $c\in(1,c_p)\cup(\frac{p}{p-3},\infty).$\footnote{Here $\tfrac p{p-3}=\infty,$ if $p=3.$} Also it is worth to note that Hakkaev, Iliev and Kirchev in \cite{hakkaev1} studied the orbital stability of a type of generalized BBM and Camassa Holm equations. The family of BBM equations that they investigated were of the form
\[u_t+2\omega u_x+3uu_x-u_{xxt}=0, \ \ \ \ \omega\in\mathbb{R}.\]
 They proved the existence solutions of the cnoidal type, but they only proved the orbital stability of this solutions in the case $\omega=0.$\\

Finally, this paper is organized as follows: In Section 2 we introduced some notations to be used throughout the whole article; in Section 3, we prove the global well-posedness result in the periodic setting; in section 4, the ill-posedness result is obtained; in Section 5, we show the existence of periodic travelling waves using the Poisson Summation theorem, in Section 6, the spectral properties needed to obtain the nonlinear stability are gotten, in Section 7, we get the stability of the  waves based on the ideas in \cite{Albert1},\cite{anguloNatali}, \cite{benjamin1} and \cite{weinstein2}, in Section 8, we present the extension of the theory, then we used it to proved the stability of cnoidal waves associated to the BBM equation.


\section{Notation and preliminaries}
The $L^2$-based Sobolev spaces of periodic functions are defined as follows (for further details see \cite{ioriolibro}). Let $\mathcal{P}=C^{\infty}_{per}$ denote the collection of all functions $f:\mathbb{R}\rightarrow \mathbb{C}$ which are $C^{\infty}$ and periodic with period $2L>0.$ The collection $\mathcal{P}'$ of all continuous linear functionals from $\mathcal{P}$ into $\mathbb{C}$ is the set of \textit{periodic distributions.} If $\Psi\in \mathcal{P}'$ then we denote the value of $\Psi$ at $\varphi$
by $\Psi(\varphi)=\langle\Psi,\varphi\rangle.$ Define the functions $\Theta_k(x)=\exp(\pi ikx/L), \ k\in \mathbb{Z},\ x\in\mathbb{R}.$ The Fourier transform of $\Psi$ is the function $\widehat{\Psi}:\mathbb{Z}\rightarrow\mathbb{C}$ defined by the formula $\widehat{\Psi}(k)=\frac 1{2L}\langle\Psi,\varphi\rangle, \ k\in\mathbb{Z}.$ So, if $\Psi$ is a periodic function with period $2L,$  we have 
\[\widehat{\Psi}(k)=\frac 1{2L}\int_{-L}^L \Psi(x)e^{-\frac{ik\pi x}{L}}dx.\]
For $s\in \mathbb{R},$ the Sobolev space of order $s,$ denoted by $H^s_{per}([-L,L])$ is the set of all $f\in \mathcal{P}'$ such that $(1+|k|^{2})^s|\widehat{f}(k)|^2\in l^2(\mathbb{Z}),$ with norm
\[||f||^2_{H^s_{per}}=2L\sum_{k=-\infty}^{\infty}(1+|k|^{2})^s|\widehat{f}(k)|^2.\]
Recall the inversion formula
\[f(x) = \sum_{k=-\infty}^{\infty}e^{\frac{ikx\pi}{L}}\widehat{f}(k),\]
and the convolution formula
\[(\widehat{f}*\widehat{g})(k) = \widehat{fg}(k),\]
where the convolution of two functions $\widehat{f}$ and $\widehat{g}$ on $\mathbb{Z}$ is formally defined by
\[(\widehat{f}*\widehat{g})(k) = \sum_{n=-\infty}^{\infty}\widehat{f}(k-n)\widehat{g}(n).\]
We also note that  $H^s_{per}$ is a Hilbert space with respect to the inner product
\[(f|g)_s = 2L\sum_{n=-\infty}^{\infty}(1+|k|^2)^{s}\widehat{f}(k)\overline{\widehat{g}(k)}\]
In the case $s=0, H^0_{per} $ is a Hilbert space that is isometrically isomorphic to $L^2([-L,L])$ and
\[(f|g)_0 = (f,g) = \int_{-L}^{L} f\overline{g} \ dx.\]
The space $H^0_{per}$ will be denoted by $L^2_{per}$ and its norm will be $\|\cdot\|_{L^2_{per}}.$
Of course $H^s_{per} \subset L^2_{per}$, for any $s \geq 0 $. Moreover, $(H^s_{per})'$, the topological dual of $H^s_{per}$, is isometrically isomorphic to $H^{-s}_{per}$ for all $s \in \mathbb{R}$. The duality is implemented concretely by the pairing
\[\langle f,g\rangle_s = 2L\sum_{k=-\infty}^{\infty}\widehat{f}(k)\overline{\widehat{g}(k)}, \ \ \ for \ \ \ f \in H^{-s}_{per}, \ \  g \in H^s_{per}. \]
Thus, if $f \in L^2_{per}$ and $g \in H^s_{per} $, with $s\geq 0,$ it follows that $\langle f,g\rangle_s = (f,g)$. Additionally, in the particular case $s=\frac 1{2}$ we will denote the pairing $\langle f,g\rangle_s$ simply by $\langle f,g\rangle.$ One of Sobolev's Lemmas in this context states that if $s>\frac{1}{2}$ and
\[C_{per} = \{f: \mathbb{R} \longrightarrow \mathbb{C} \ | \ f \ \ \text{is continuous and periodic with period} \ \  2L \},\]
then $H^{s}_{per}\hookrightarrow C_{per}$.\\

If $Y$ is a Banach space like $H^{s}_{per}$, and $T > 0$, then $C([0,T];Y)$ is the space of continuous mappings from $[0,T]$ to $Y$ and, for $k \geq 0,\ C^k([0,T];Y)$ is the subspace of mappings $t\mapsto u(t)$ such that $\partial^j_{t}u \in C([0,T];Y)$ for $0\leq j\leq k$, where the derivative is taken in the sense of vector-valued distributions. This space carries the standard norm
\[\|u\|_{C^k([0,T];Y)} = \sum_{j=0}^{k} \max_{0\leq t\leq T}\|\partial^j_{t}u(t)\|_Y.\]
$c_0$ denotes numeric positive constants which can change form line to line. For any positive elements  $A$ and $B$ the notation $A\apprle B$ (respectively, $A\apprge B$) means that there exist a positive constant  $c_0$ such that $A\leq c_0B$ (respectively, $A\geq c_0B$).
The notation  $A\sim B$ means that $A\apprle B\apprle A.$ Finally $\mu(A)$ denotes the Lebesgue measure of the set $A.$\\

Next, we present the Poisson Summation theorem. It will be used in Section 5 and 8 to find the periodic travelling wave solutions for the rBO and the BBM equation, respectively.
\begin{theo}\label{theoSumPois}
Let ${\widehat{f}}^{\mathbb{R}}(\xi)=\int_{-\infty}^{\infty}f(x)e^{-ix\xi}dx$ and $f(x)=\int_{-\infty}^{\infty}\widehat{f}^{\mathbb{R}}(\xi)e^{ix\xi}d\xi$ satisfy
\[|f(x)|\leq\frac{A}{(1+|x|)^{1+\delta}} \ \ \ \ \text{and} \ \ \ |\widehat{f}^{\mathbb{R}}(\xi)|\leq\frac{A}{(1+|\xi|)^{1+\delta}},\]
where $A>0$ and $\delta>0$ (then $f$ and $\widehat{f}$ can be assumed continuous functions). Thus, for $L>0$
\[\sum_{n=-\infty}^{\infty}f\left(x+2Ln\right)=\frac{1}{2L}\sum_{n=-\infty}^{\infty}\widehat{f}^{\mathbb{R}}\left(\frac{n}
{2L}\right)e^{\frac{\pi inx}{L}}.\]
The two series above converge absolutely.
\end{theo}
\proof
See for example \cite{steinWiess}.
\endproof

\section{Well-posedness results}
In this section we study the well-posedness for the periodic initial value problem (\ref{probperiod}). For simplicity in this section we wiil consider $L=\pi.$\\
First, rewrite (\ref{ecuabasica}) as
\[(1+\mathcal{H}\partial_x)u_t=-(u+\tfrac 1{2}u^2)_x,\]
since  $\mathcal{H}\partial_x \geq 0$, formally we have
\[u_t=-\partial_x\left(1+\mathcal{H}\partial_x\right)^{-1}\left(u+\tfrac 1{2}u^2\right)=K(u+\tfrac 1{2}u^2),\]
where $K$ is given explicitly by its Fourier transform as
\begin{equation}\label{formadeK}
 \widehat{Ku}(n)=\frac{-in}{1+|n|}\widehat{u}(n).
\end{equation}
Integrating and using the initial condition we get
\[u(x,t)=u_0(x)+\int_{0}^{t} K(u+\tfrac 1{2}u^2)(x,\tau)d\tau,\]
for all $x\in \mathbb{R}$ and $t>0.$\\
Let $s>\frac 1{2}$ and $T>0,$  define $A$ on $X= C([0,T];H_{per}^s),$ by 
\[Au(t)=u_0+\int_{0}^{t} K(u+\tfrac 1{2}u^2)(\cdot ,\tau)d\tau.\] 
Using (\ref{formadeK}) and the fact that  $H_{per}^s$ with $s>\frac 1{2}$ is Banach algebra we obtain
\begin{align}\label{desioperador1}
\|Au(t)\|_{ H_{per}^s} \leq \|u_0\|_X+T\left(\|u\|_X+\tfrac{c_0}{2}\|u\|^2_X\right).
\end{align}
On the other hand, if $u,w\in X,$ then for  $0\leq t\leq T$ we get
\begin{align}\label{desioperador2}
\|Au(t)-Aw(t)\|_{ H_{per}^s} &\leq T\|u-w\|_X \left[1+\frac {c_0}{2}(\|u\|_X+\|w\|_X)\right].
\end{align}
Now, define $M=\{u\in X:\|u\|_X\leq R\}.$  If $u,w\in M$ from (\ref{desioperador1}) and (\ref{desioperador2}) we conclude 
\[\|Au\|_X \leq \|u_0\|_X+T(R+\frac {c_0}{2}R^2)\]  and
\[\|Au-Aw\|_X\leq T\|u-w\|_X (1+c_0R).\]
Taking  $R=2\|u_0\|_X$ and \ $T=\frac 1{2}(1+c_0R)^{-1},$ it follows that  $A:M\longrightarrow M$ is a contraction, hence there exists a unique $u\in M$ such that $Au(t)=u(t)$ for all $t\in [0,T].$ Collecting this information we have the next theorem.
\begin{theo}
Suppose $s> \frac{1}{2},$ then for all $u_0\in H_{per}^s([-\pi,\pi])$ there exists $T=T(\|u_0\|_{H_{per}^s})>0$  and a unique solution of (\ref{probperiod}) on the interval $[-T,T]$, such that $u \in C ([-T,T];H_{per}^s)$. Furthermore, for all $T'<T$ there exists a neighborhood $V$ of $u_0$ in $H_{per}^s([-\pi,\pi])$ such that 
\[\mathbb{F}:V \longrightarrow C([-T',T'];H_{per}^s([-\pi,\pi])), \ \tilde{u}_0 \rightarrow \tilde{u}(t),\]
is Lipschitz.
\end{theo}
Now, we study the global well-posedness. The following lemmas are useful to derive an \textit{a priori} estimate.
\begin{lemm}\label{lema1} Suppose $u_0\in H_{per}^s,$ where  $s>\frac 1{2}$ and let $u$ be the corresponding solutions of  (\ref{probperiod}) on the interval $[0,T]$. Then for all $t\in [0,T]$
\[\|u(t)\|_{ H_{per}^{\frac{1}{2}}}=\|u_0\|_{ H_{per}^{\frac{1}{2}}}.\]
\end{lemm}
\begin{lemm} \label{lema2} If $ s_0\leq s \leq s_1$, with $s=\theta s_0+(1-\theta)s_1,$ and $0\leq \theta \leq 1$  , then there exists a constant $C>0,$ such that 
\begin{equation} \label{estimativacom1}
\|J^sf\|_{L_{per}^2}\leq C\|J^{s_0}f\|^{\theta}_{L_{per}^2}\|J^{s_1}f\|^{1-\theta}_{L_{per}^2, } 
\end{equation}
where $J^s=(1+\Delta)^{\frac s{2}}$ ( The Bessel potencial of order $s$).
\end{lemm}
\begin{lemm} \label{lema3} Let  $s\geq 1$, then there exists a constant $C=C(s)>0$  such that 
\begin{equation} \label{estiamtivacom2}
\|J^s(fg)-fJ^sg\|_{L_{per}^2}\leq C\left\{\|\partial_xf\|_A\|g\|_{ H_{per}^{s-1}}+ \|\partial_xf\|_{ H_{per}^{s-1}}\|g\|_A\right\} 
\end{equation}
where $||g||_A=\sum_{k=-\infty}^{\infty}|\widehat{g}(k)|.$
\end{lemm}
The proof of Lemma \ref{lema1} is immediate from the conservation law $F$ in (\ref{prisegley}) (see for instance \cite{bona1}), the proof of Lemma \ref{lema2} can be found in \cite{bergh1} and Lemma \ref{lema3} is a commutator estimative, see Lemma B.3 in \cite{ioriolibro}.  Next, we establish the promised estimate.
\begin{theo}
Let $s\geq 3/2.$ If $u\in C([0,T],H_{per}^s)$ is a solution of (\ref{probperiod}) in $[0,T]\times \mathbb{R},$ with initial datum $u_0,$ then there exist constants $C=C(\|u_0\|_{ H_{per}^{\frac{3}{2}}})>0$ and $C_s=C_s(\|u_0\|_{ H_{per}^s})>0$ such that
\begin{equation} \label{desiexp}
\sup_{t\in [0,T]}\|u(t)\|_{ H_{per}^s}\leq C_s e^{CT}
\end{equation}
\end{theo}
\proof
First we will establish the result in the case $s=3/2,$ for this consider a sequence $\{u^n\}_{n=0}^{\infty}$ in $C^1([0,T],H_{per}^{\infty}),$ converging to $u$ in  $C^1([0,T],H_{per}^{\frac{3}{2}})$ and define  $F(v)=v_t+v_x+vv_x+\mathcal{H}v_{xt} .$
Then $F(u^n)\rightarrow 0$ in $C([0,T],H_{per}^{\frac{1}{2}})$ and $\{\partial_x^2u^n\}_{n=0}^{\infty}$ is bounded in $ C([0,T],H_{per}^{-\frac 1{2}}).$ In fact,
\begin{align*}
\|\partial_x^2u^n(t)\|^2_{ H_{per}^{-\frac 1{2}}}&\leq 2\pi \sum_{k=-\infty}^{\infty}(1+|k|^2)^{-\frac 1{2}}|k|^4|\widehat{u^n}(k)|^2\leq 2\pi \sum_{k=-\infty}^{\infty}(1+|k|^2)^{\frac {3}{2}}|\widehat{u^n}(k)|^2\\&=\|u^n(t)\|^2_{ H_{per}^{\frac {3}{2}}}
\end{align*}
and
\begin{align*}
\|F(u^n)\|_{H^{\frac{1}{2}}_{per}} &\leq \|u_{t}^n - u_{t}\|_{H_{per}^{\frac{1}{2}}} + \|u_{x}^n - u_{x}\|_{H_{per}^{\frac{1}{2}}} + \|\mathcal{H}(u_{xt}^n - u_{xt})\|_{H_{per}^{\frac{1}{2}}} + \|u^nu_{x}^n - uu_{x}\|_{H_{per}^{\frac{1}{2}}}\\
& \leq 2\|u_{t}^n - u_{t}\|_{H_{per}^{\frac{3}{2}}} + \|u^n - u\|_{H_{per}^{\frac{3}{2}}} + \|u_{x}^n(u^n - u)\|_{H_{per}^{\frac{1}{2}}}+ \|u(u_{x}^n - u_{x})\|_{H_{per}^{\frac{1}{2}}}\\
&\leq 2\|u_{t}^n - u_{t}\|_{H_{per}^{\frac{3}{2}}} + \|u^n - u\|_{H_{per}^{\frac{3}{2}}} + c_0\|u_{x}^n \|_{H_{per}^{\frac{1}{2}}} \|u^n - u\|_{H_{per}^{\frac{3}{2}}}\\&\ + c_0\|u\|_{H_{per}^{\frac{3}{2}}}\|u_{x}^n- u_{x}\|_{H_{per}^{\frac{1}{2}}}\leq 2\|u_{t}^n - u_{t}\|_{H_{per}^{\frac{3}{2}}} + \|u^n - u\|_{H_{per}^{\frac{3}{2}}} \\&\ +c_0\|u^n \|_{H_{per}^{\frac{3}{2}}} \|u^n - u\|_{H_{per}^{\frac{3}{2}}} + c_0\|u\|_{H_{per}^{\frac{3}{2}}}\|u^n - u\|_{H_{per}^{\frac{3}{2}}},
\end{align*}
where we used the fact that $\|uv\|_{H_{per}^{\frac{1}{2}}}\leq c_0\|u\|_{H_{per}^{\frac{1}{2}}} \|v\|_{H^{s}_{per}} $, if $ u \in {H_{per}^{\frac{1}{2}}} $ and $v \in H^s_{per}$ with $s>1$.
Hence, $\langle F(u^n),\partial_x^2u^n\rangle \rightarrow 0,$ when $n\rightarrow \infty,$ in $C([0,T]).$\\

On the other hand, since $u^n$ is smooth, we have that
\begin{align} \label{desifconsegunderi}
2\langle F(u^n),\partial_x^2u^n\rangle=-\frac d{dt}\|u_x^n(\cdot,t)\|_{ H_{per}^{\frac 1{2}}} -\int_{-\pi}^{\pi}(u_x^n)^3dx.
\end{align}
The Sobolev inequality and  (\ref{estimativacom1}) imply that there exists a positive constant $c_0$ such that 
\begin{align} \label{desiderifalcubo}
\int_{-\pi}^{\pi}(\partial_xf)^3dx \notag &\leq \|\partial_xf\|^3_{L_{per}^3}\leq c_0\|\partial_xf\|^3_{H_{per}^{\frac 1{6}}}\leq c_0\|f\|^3_{H_{per}^{\frac 7{6}}}\leq c_0\|f\|^{3\frac 1{3}}_{H_{per}^{\frac 1{2}}}\|f\|^{3(1-\frac 1{3})}_{H_{per}^{\frac 3{2}}}\\&=c_0 \|f\|_{ H_{per}^{\frac 1{2}}} \|f\|^2_{ H_{per}^{\frac 3{2}}}
\end{align}
Integrating (\ref{desifconsegunderi}) form  $0$ to $t$ and using (\ref{desiderifalcubo}) we get
\begin{align*}
2\int_{0}^{t}\langle F(u^n),\partial_x^2u^n\rangle d\tau \leq &\ \|u_x^n(\cdot,0)\|^2_{ H_{per}^{\frac 1{2}}}-\|u_x^n(\cdot,t)\|^2_{ H_{per}^{\frac 1{2}}}\\& +c_0\int_{0}^{t}\|u^n(\cdot,\tau)\|_{ H_{per}^{\frac 1{2}}}\|u^n(\cdot,\tau)\|^2_{ H_{per}^{\frac 3{2}}}d\tau.
\end{align*}
Using the last inequality we obtain,
\begin{align*}
\|u^n(\cdot,t)\|^2_{ H_{per}^{\frac 3{2}}}&\leq \|u^n(\cdot,t)\|^2_{ H_{per}^{\frac 1{2}}}+\|u_x^n(\cdot,t)\|^2_{ H_{per}^{\frac 1{2}}}\leq \|u^n(\cdot,t)\|^2_{ H_{per}^{\frac 1{2}}}+ \|u_x^n(\cdot,0)\|^2_{ H_{per}^{\frac 1{2}}}\\&+c_0\int_{0}^{t}\|u^n(\cdot,\tau)\|_{ H_{per}^{\frac 1{2}}}\|u^n(\cdot,\tau)\|^2_{ H_{per}^{\frac 3{2}}}d\tau -2\int_{0}^{t}\langle F(u^n),\partial_x^2u^n\rangle d\tau.
\end{align*}
Taking the limit when $n\rightarrow \infty$ and using Lemma \ref{lema1}, yields
\[\|u(\cdot,t)\|^2_{ H_{per}^{\frac 3{2}}}\leq 2\|u_0\|^2_{ H_{per}^{\frac 3{2}}}+c_0\|u_0\|_{ H_{per}^{\frac 1{2}}}\int_{0}^{t}\|u(\cdot,\tau)\|^2_{ H_{per}^{\frac 3{2}}}d\tau.\]
Finally by the Gronwall inequality we have that
 \[\|u(\cdot,t)\|^2_{ H_{per}^{\frac 3{2}}}\leq  2\|u_0\|^2_{ H_{per}^{\frac 3{2}}}\exp\left(c_1t\|u_0\|_{ H_{per}^{\frac 1{2}}}\right)\]
for all $t\in [0,T],$ so we have the result in the case $s=3/2.$\\

For the general case suppose that (\ref{desiexp}) is true for $r_0\geq 3/2,$ and let $s=r_0+\alpha,$ with $0<\alpha \leq 1/2.$ Additionally, with the only goal of making the reading easier, let us define $r=s-1/2.$ We will prove that the inequality holds for $s.$ Consider $\{u^n\}_{n=0}^{\infty}$ a sequence in $C^1([0,T],H_{per}^{\infty})$ such that
$u^n\rightarrow u,$ in $C^1([0,T],H_{per}^s),$ then $F(u^n)\rightarrow 0,$  in $C([0,T],H_{per}^{s-1})$ and $\{J^ru^n\}_{n=0}^{\infty}$ is bounded in $C([0,T],H_{per}^{s-r}).$ Therefore, $\langle J^ru^n,J^rF(u^n)\rangle \longrightarrow 0$ in $C([0,T]),$ when $n\rightarrow \infty.$\\

Let us define $v:=u^n,$ it easy to see that
\[2\langle J^{r}v,J^rF(v)\rangle = \frac{d}{dt} \langle J^{r}v,J^r{v}+J^{r}\mathcal{H}v_{x}\rangle + 2 \langle J^{r}v,J^r(vv_{x})\rangle. \]
Then, integrating from  $0$ to $t$ the last equality, using the fact that
\begin{align} \label{igualconj}
\langle J^{r}v,J^r(vv_{x}) \rangle &= -\frac{1}{2} (J^{r}v,v_{x}J^{r}v) + \langle J^{r}v,J^r(vv_{x})-vJ^rv_{x} \rangle
\end{align}
and since $s = r +1/2,$ one has that 
\begin{align*}
\|v(\cdot,t)\|^2_{ H_{per}^s} \leq &\ 2L\sum (1+|k|^2)^r(1+|k|) |\widehat{v}(k)|^2 
\\=&\ \langle J^{r}v,J^rv \rangle + \langle J^{r}v,J^r\mathcal{H}v_{x}\rangle
\\=&\ \langle J^{r}v(0),J^rv(0)+ J^r\mathcal{H}v_{x}(0) \rangle + 2 \int_{0}^{t}\langle J^{r}v,J^rF(v)\rangle d\tau\\&-2\int_{0}^{t}\langle J^{r}v,J^r(vv_{x})\rangle d\tau.
\end{align*}
Taking the limit when $n \longrightarrow \infty$, we obtain
\begin{equation} \label{desiconjyu}
\|u(\cdot,t)\|^2_{ H_{per}^s} \leq c_s \|u_0\|^2_{ H_{per}^r} - 2 \lim_{n\rightarrow \infty}\int_{0}^{t} \langle J^{r}v,J^r(vv_{x})\rangle d\tau.
\end{equation}
Using (\ref{igualconj}) and the Lemma \ref{lema3}
\begin{align*}
2 |\langle J^{r}v,J^r(vv_{x})\rangle| &\leq |( J^{r}v,v_{x}J^{r}v))|+ 2\|J^{r}v\|_{L_{per}^2}\|J^{r}(vv_{x})-vJ^{r}v_x\|_{L_{per}^2}  
\\& \leq \|J^{r}v\|^2_{L_{per}^2}\|v_{x}\|_{\infty} + 2\|J^{r}v\|_{L_{per}^2}2c_0\|v_{x}\|_A \|v_{x}\|_{H_{per}^{r-1}}
\\& \leq c_0\|v\|^2_{H_{per}^r}\|v\|_{H_{per}^s} + 4c_0\|v\|_{H_{per}^r}\|v_x\|_{H_{per}^{s-1}}\|v\|_{H_{per}^{r}} 
\\& \leq c_0\|v\|^2_{H_{per}^{r}}\|v\|_{H_{per}^{s}}.
\end{align*}
Going back to (\ref{desiconjyu}) it is obtained,
\[\|u(\cdot,t)\|^2_{ H_{per}^s} \leq c_s\|u_0\|^2_{ H_{per}^r} + c_0\int_{0}^{t}\|u(\cdot,\tau)\|^2_{ H_{per}^r}\|u(\cdot,\tau)\|_{ H_{per}^s}d\tau.\]
Now, since $s>r$ and  $r_0\geq r,$ then
\begin{align*}\|u(\cdot,t)\|^2_{ H_{per}^s} &\leq c_s \|u_0\|^2_{ H_{per}^s} + c_0\int_{0}^{t} \|u(\cdot,\tau)\|_{ H_{per}^r} \|u(\cdot,\tau)\|_{ H_{per}^r}\|u(\cdot,\tau)\|_{ H_{per}^s}d\tau \\&\leq c_s \|u_0\|^2_{ H_{per}^s} + c_0\int_{0}^{t} \|u(\cdot,\tau)\|_{ H_{per}^{r_0}}\|u(\cdot,\tau)\|^2_{ H_{per}^s} d\tau.
\end{align*}
So, using the inductive hypothesis
\[\|u(\cdot,t)\|^2_{ H_{per}^s} \leq  c_s\|u_0\|^2_{ H_{per}^s} + C_{r_{0}}e^{c_0T} \int_{0}^{t}\|u(\cdot,\tau)\|^2_{ H_{per}^s}d\tau,\]\
we obtain the desired result by the Gronwall inequality.
\endproof
\begin{coro}
The periodic initial value problem (\ref{probperiod}) is globally well-posed in $H_{per}^s([-\pi,\pi])$ for $s\geq 3/2.$ 
\end{coro}


\section {Ill-posedness results}
Since arguments of scaling cannot be applied to the rBO equation to obtain a notion of criticality, it is not clear which  is the smaller indice for we can have well-posedness for this equation in the Sobolev space $H^s_{per}$ (or $H^s(\mathbb{R})).$ In this section we will partially answer this question, we will show that the flow map data-solution for the Cauchy problem associated to the rBO equation is not smooth (it is not $C^2$) at the origin for initial data in $H^s_{per}$ (or $H^s(\mathbb{R})$), with  $s<0.$ Therefore we cannot apply the Contraction Principle to solve the integral equation (\ref{eqInteAsso}), see for instance \cite{LinaresPonce1}.\\

First we will analyze the problem on the periodic setting. For simplicity we consider functions of period $2\pi.$ We know that the linear problem associated to (\ref{ecuabasica}) with initial datum $\phi$ has solution 
\[u(x,t)=S(t)\phi(x),\]
where $S$ satifies
\[\widehat{S(t)\phi}(n)=e^{-\tfrac{in}{1+|n|}t}\widehat{\phi}(n).\]
Now, if $u$ is solution of (\ref{ecuabasica}), then by the Duhamel principle we have that
\begin{equation}\label{eqInteAsso}
 u(x,t)=S(t)\phi(x)-\int_{0}^{t}S(t-\tau)\Lambda[u(x,\tau)u_{x}(x,\tau)]d\tau,
\end{equation}
where $\widehat{\Lambda u}(n)=(1+|n|)^{-1}\widehat{u}(n).$\\

Next, we prove the principal result of this section.
\begin{theo}\label{t1illpos}
Let $s<0$ and  $T$ a positive number. Then there does not exist a space $X_T$ continuously embedded in $C([-T,T];H_{per}^{s})$ such that there exist $c_0>0$ satisfying
\begin{equation} \label{teo1a}
\|S(t)\phi\|_{X_T}\leq c_0\|\phi\|_{H_{per}^{s}}, \,\,\,\ \forall\phi\in H^s_{per} 
\end{equation}
and
\begin{equation} \label{teo1b}
\left|\left|\int_{0}^{t}S(t-\tau)\Lambda [u_x(\tau)u(\tau)]d\tau\right|\right|_{{X}_{T}} \leq c_0\|u\|^2_{{X}_{T}}  \,\,\,\ \forall u \in {X}_{T}.
\end{equation}
\end{theo}
\proof
Suppose that there exists a space $X_T$ continuously embedded in $C([-T,T];H_{per}^{s})$ such that (\ref{teo1a}) and (\ref{teo1b}) hold. Consider $\phi \in H_{per}^{s}$ and define $u := S(t)\phi,$ by the assumption we have that  $u\in X_T$ and 
\[\left|\left|\int_{0}^{t}S(t-\tau)\Lambda[S(t)\phi (S(t)\phi)_x]d\tau\right|\right|_{{X}_{T}} \leq c_0\|S(t)\phi\|^2_{{X}_{T}}\apprle\|\phi\|^2_{H_{per}^{s}}. \]
Now, since $X_T \hookrightarrow C([-T,T];H_{per}^{s}),$ we get 
\begin{equation} \label{equa3}
\left|\left|\int_{0}^{t}S(t-\tau)\Lambda[S(t)\phi(S(t)\phi)_x]d\tau\right|\right|_{H_{per}^{s}}\leq c_0\|\phi\|^2_{H_{per}^{s}}.
\end{equation}
We will prove that (\ref{equa3}) does not hold choosing $\phi$ appropriate. For this, consider
\[\phi(x):= N^{-s}\cos(Nx),\ \ \  \text{with} \ \ N\in \mathbb{N},\ \  N\gg 1.\]
First, it easy to see that   
\begin{align*}
S(t)\phi(x) &=N^{-s}\cos\left(Nx-\frac{N}{1+N}t\right).
\end{align*}
Then, 
\begin{align*}
\psi(x,t)&:= \int_{0}^{t}S(t-\tau)\Lambda[S(t)\phi(x) (S(t)\phi(x))_x]d\tau\\
&=-\frac{1}{2}N^{-2s+1}\int_{0}^{t}S(t-\tau)\Lambda\left[\sin\left(2Nx-\tfrac{2N}{1+N}\tau
\right)\right]d\tau.
\end{align*}
Now, using the specific form of $\Lambda$ we obtain that
\begin{align*}
\int_{0}^{t}S(t-\tau)\Lambda\left[\sin\left(2Nx-\tfrac{2N}{1+N}\tau
\right)\right]d\tau&= -\frac{1}{2(1+2N)\gamma_N}\left[ e^{i\left( 2Nx-\frac{2N}{1+2N}t\right)} -e^{i\left( 2Nx-\frac{2N}{1+N}t\right)}\right] \\ 
& +\frac{1}{2(1+2N)\gamma_N}\left[e^{-i\left( 2Nx-\frac{2N}{1+N}t\right)}- e^{-i\left( 2Nx-\frac{2N}{1+2N}t\right)}\right],
\end{align*}
where $\gamma_N = \frac{2N^2}{(1+N)(1+2N)}.$
Therefore 
\[\psi(x,t) = \frac{1}{2}N^{-2s+1}\frac{1}{\gamma_N(1+2N)}\left[\cos\left(2Nx -\frac{2N}{1+2N}	t\right)-\cos\left(2Nx -\frac{2N}{1+N}	t\right)\right].\]

Now, it is easy to see that
\[\|\psi(\cdot,t)\|_{H_{per}^s}^2 \sim N^{-4s}\left| e^{-i\frac{2N}{1+2N}t} - e^{-i\frac{2N}{1+N}t}\right|^{2}(1+4N^2)^s.\]
But, $\left|e^{-i\frac{2N}{1+2N}t}-e^{-i\frac{2N}{1+N}t}\right|^{2}= 2 - 2\cos\left(\frac{2N}{1+N}t - \frac{2N}{1+2N}t\right).$
Hence,
\[\|\psi(\cdot,t)\|_{H_{per}^2} \sim N^{-s}\left( 1 - \cos \left(\gamma_Nt\right)\right)^{\frac{1}{2}}.\]
Note that $\|\phi\|^2_{H^s_{per}}\sim 1,$ then for all $t\in (0,T)$ we have
\[\frac{\|\psi(\cdot,t)\|_{H_{per}^s}}{\|\phi\|^2_{H^s_{per}}}\sim N^{-s}\left( 1 - \cos \left(\gamma_Nt\right)\right)^{\frac{1}{2}}.\]
Without loss of generality suppose that $0<T<2\pi$ and consider $s<0$ fixed, then since $\gamma_N \rightarrow 1^{-},$ as $N \rightarrow +\infty,$ we obtain for all $0<t<T$ that
\[\frac{\|\psi(\cdot,t)\|_{H_{per}^s}}{\|\phi\|^2_{H^s_{per}}} \longrightarrow +\infty,\] 
as $N \rightarrow +\infty,$ but this contradict (\ref{equa3}), which finishes the proof of the theorem.
\endproof

As a consequence of the last theorem we get the next result.
\begin{theo}\label{t2illpos}
Fix $s<0.$ There does not exist a  $T>0$ such that (\ref{ecuabasica}) admits a unique local solution define on the interval $[-T,T]$ and such that the flow map data-solution 
\[\phi \longmapsto u(t), \,\,\,\ t \in [-T,T] \]
for (\ref{ecuabasica}) is $C^2$ differentiable at zero from $H_{per}^{s}$ to $H_{per}^{s}.$
\end{theo}
\proof
Consider the Cauchy problem
\begin{equation} \label{equacongama}
\left \{ 
\begin{aligned}
u_t + u_x + uu_x+\mathcal{H}u_{xt}&=0, \\
u(0,x)&=\phi_{\gamma}(x), \ \ \ \ \ 0<\gamma\ll 1\\
\end{aligned} \right.
\end{equation}\\
where $\phi_{\gamma}(x)=\gamma\phi(x).$ Suppose that $u(\gamma,t,x)$ is a local solution of  (\ref{ecuabasica}) and that the flow map data-solution is  $C^2$ at the origin from $H_{per}^s$ to $H_{per}^s$. Then
\[\frac{\partial u}{\partial \gamma}(\gamma,t,x)|_{\gamma=0} = S(t)\phi(x)\]
and
\[\frac{\partial^2 u}{\partial \gamma}(\gamma,t,x)|_{\gamma=0} = -2\int_{0}^{t}S(t-\tau)\Lambda\left[(S(\tau)\phi)(S(\tau)\phi)_x\right]d\tau.\]
Using again the assumption that the flow map is  $C^2$ at zero, we have
\[\left|\left|\int_{0}^{t}S(t-\tau)\Lambda[(S(\tau)\phi)(S(\tau)\phi)_x]d\tau\right|\right|_{H_{per}^s} \leq c_0\|\phi\|_{H_{per}^s}^2.\]
But the last estimative is the same in (\ref{equa3}), which has been shown to fail in the last theorem. This finishes the proof.
\endproof

Now we will establish the result in the non-periodic setting. Recall that in this case we have  
\[S(t)\phi(x)=\int_{\mathbb{R}}\widehat{\phi}(\xi)e^{i\left(\xi x-{\tfrac{\xi}{1+|\xi|}t}\right)}d\xi\]
and $\widehat{\Lambda u}(\xi)=(1+|\xi|)^{-1}\widehat{u}(\xi).$ Let us start with a lemma.
\begin{lemm}
\begin{align*}
\int_{0}^{t}S(t-\tau)\Lambda[(S(\tau)\phi)&(S(\tau)\phi)_x]d\tau =\\&c_0\int_{\mathbb{R}^2}e^{i(\xi x-p(\xi)t)}\frac{\xi}{1+|\xi|}\widehat{\phi}(\eta)\widehat{\phi}(\xi-\eta)\frac{e^{-it\chi(\xi,\eta)}-1}{\chi(\xi,\eta)}d\eta d\xi
\end{align*}
where $p(\xi)=\frac{\xi}{1+|\xi|}$ and $\chi(\xi,\eta)=p(\eta)+p(\xi-\eta)-p(\xi).$
\end{lemm}
\proof 
The proof is very similar to Lemma 1 in \cite{MoliSautTzve1}, we will do it here just for sake of completeness. Using the inverse Fourier transform is easy to see that
\begin{align*}
&\int_{0}^{t}S(t-\tau)\Lambda[(S(\tau)\phi)(S(\tau)\phi)_x]d\tau\\
&=\int_{0}^{t}\int_{\mathbb{R}}e^{i(\xi x-p(\xi)t)}e^{\tau p(\xi)}\left[\frac{i\xi}{1+|\xi|}(S(\tau)\phi)\widehat * (S(\tau)\phi)\widehat (\xi)\right] d\xi d\tau\\
&= i \int_{0}^{t}\int_{\mathbb{R}}e^{i(\xi x-p(\xi)t)}e^{\tau p(\xi)}\frac{\xi}{1+|\xi|}[(e^{-i\tau p(\cdot)}\widehat{\phi}(\cdot))*(e^{-i\tau p(\cdot)}\widehat{\phi}(\cdot))] (\xi)d\xi d\tau\\
&=i\int_{\mathbb{R}}e^{i(\xi x-p(\xi)t)}\frac{\xi}{1+|\xi|}\widehat{\phi}(\eta)\widehat{\phi}(\xi-\eta)\int_{0}^{t}e^{i\tau[p(\eta)+p(\xi-n)-p(\xi)]} d\tau d\eta d\xi\\
&=i\int_{\mathbb{R}^2}e^{i(\xi x-p(\xi)t)}\frac{\xi}{1+|\xi|}\widehat{\phi}(\eta)\widehat{\phi}(\xi-\eta)\frac{e^{-i\tau[p(\eta)+p(\xi-n)-p(\xi)]}-1}{p(\eta)+p(\xi-n)-p(\xi)} d\eta d\xi, 
\end{align*}
which finishes the proof of the lemma.
\endproof

Next, define
\[\varphi(x,t):= \int_{0}^{t}S(t-\tau)\Lambda[(S(\tau)\phi)(S(\tau)\phi)_x]d\tau,\]
then using the last lemma we have that
\begin{equation}\label{ma1}
\widehat{\varphi}(\xi,t)=c_0\frac{\xi}{1+|\xi|}e^{-ip(\xi)t}\int_{\mathbb{R}}
\widehat{\phi}(\eta)\widehat{\phi}(\xi-\eta)\frac{e^{-it\chi(\xi,\eta)}-1}{\chi(\xi,\eta)}d\eta.
\end{equation}
In this case we consider
\[\widehat{\phi}(\xi):=N^{-s}\chi_{[N,N+1]}(\xi), \ \ \ \text{with}\ \  N \in \mathbb{N},\ \  N\gg 1\]
where $\chi_A$ denotes the caracteristic function of $A.$ Note that $\|\phi\|_{H^s(\mathbb{R})} \sim 1$. Additionally, using (\ref{ma1}) we get
\[\widehat{\varphi}(\xi,t)=c_0\frac{\xi}{1+|\xi|}e^{-p(\xi)t}N^{-2s}\int_{\Omega_\xi}
\frac{e^{-it\chi(\xi,\eta)}-1}{\chi(\xi,\eta)}d\eta,\]
with $\Omega_\xi =\{\eta: \eta\in\text{supp} \ \widehat{\phi} \ \ \ \text{and} \ \ \ \xi-\eta\in \text{supp}\ \widehat{\phi}\}.$ Now, since $s<0$, we can chose  $\epsilon>0$ such that $-s-\epsilon>0$. Consider $t= N^{-\epsilon}$ and note that for $\xi \in \left(2N+\frac{1}{2},2N+1\right)$ we have $\mu(\Omega_\xi)\gtrsim 1$.\\
It is easy to see that
\[\chi(\xi,\eta) =\frac{\eta(\xi-\eta)(2+\xi)}{(1+\eta)(1+\xi-\eta)(1+\xi)} \leq 3 \ \ \ \  \forall \   \eta, \xi-\eta \in [N,N+1].\]
Then for $N$ big enough we arrive at
\begin{align*} 
\|\varphi(\cdot,t)\|^2_{H^s(\mathbb{R})} &\gtrsim \int_{2N+\frac{1}{2}}^{2N+1}(1+|\xi|^2)^{s}N^{-4s}\frac{|\xi|^2}{(1+|\xi|)^2}|t|^2\left|\int_{\Omega_\xi}
\frac{e^{-it\chi(\xi,\eta)}-1}{t\chi(\xi,\eta)}d\eta\right|^2d\xi \\  &\gtrsim \int_{2N+\frac{1}{2}}^{2N+1}(1+|\xi|^2)^{s}N^{-4s}\frac{|\xi|^2}{(1+|\xi|)^2}|t|^2\left|\text{Im}\int_{\Omega_\xi}
\frac{e^{-it\chi(\xi,\eta)}-1}{t\chi(\xi,\eta)}d\eta\right|^2d\xi
\\&=\int_{2N+\frac{1}{2}}^{2N+1}(1+|\xi|^2)^{s}N^{-4s}\frac{|\xi|^2}{(1+|\xi|)^2}|t|^2\left|\int_{\Omega_\xi}
\frac{\sin(t\chi(\xi,\eta))}{t\chi(\xi,\eta)}d\eta\right|^2d\xi \\&\gtrsim N^{-4s}N^{2s}t^2.
\end{align*}
Hence
\[1\sim\|\phi\|_{H^s(\mathbb{R})}\gtrsim\|\varphi(\cdot,t)\|_{H^s(\mathbb{R})} \gtrsim  N^{-s-\epsilon},\]
which is a contradiction for $N\gg 1.$ This completes the proof in the non-periodic case.\\


\section{Periodic travelling wave solutions}
The aim of this section is to show the existence of a smooth curve of periodic travelling wave solutions for 
(\ref{ecuaordina}), via the Poisson Summation theorem. In fact, consider the following equation
\[\omega\mathcal{H}\varphi_{w}'+(w -1)\varphi_{w}-\frac{1}{2}\varphi_{w}^2 = 0.\]
This equation determines solitary travelling wave solutions to the rBO equation on $\mathbb{R}$ in the form
\[\varphi_w(x)=\frac{4(w-1)}{1+\left(\frac{w-1}{w}x \right)^2} \ , \ \ \ \ \ \  w>1.\]  
Its Fourier transform is given by
\[\widehat{\varphi}^{\mathbb{R}}_w(\xi)=4\pi we^{-|\frac{w}{w-1}\xi|}.\]
Therefore, by the Poisson Summation theorem (see Theorem \ref{theoSumPois} in Section 2), we get the following periodic function,
\begin{align}\label{ecuapoisson}
\notag \psi_{w}(x) &=\sum_{n=-\infty}^{\infty}\varphi_w(x+2Ln)=\frac{2\pi w}{L}\sum_{n=-\infty}^{\infty} e^{-\frac{w|n|}{2(w-1)L}} e^{\frac{\pi inx}{L}}\\ 
\notag &=\frac{2\pi w}{L}\sum_{n=0}^{\infty}\epsilon_n\ e^{-\frac{wn}{2(w-1)L}}\cos\left(\frac{n\pi x}{L}\right)\\ 
&= \frac{2\pi w}{L}\left(\frac{\sinh\left(\frac{w}{2(w-1)L}\right)}{\cosh\left(\frac{w}{2(w-1)L}\right)-\cos \left(\frac{\pi}{L} x\right)}\right),
\end{align}
where in the last identity we used the Fourier expansion $1.89$ in  \cite{oberhettinger11} (see also \cite{anguloNatali}) and 
\[\epsilon_n=\left \{ 
\begin{aligned}
1 &, \ \ \ \text{if} \ \ n=0\\
2 &, \ \ \ \text{if} \ \ n=1,2,3, . . . \\
\end{aligned} \right.\]
Now, consider $\phi_c,$ with $c \neq 1$, a smooth periodic solution of the equation (\ref{ecuaordina}). Then $\phi_c$ and $\phi_c^2$ can be  expressed as Fourier series
\begin{equation} \label{ecaexpansao}
 \phi_c(x)=\sum_{n=-\infty}^{\infty}a_{n}e^{\frac{i\pi n x}{L}}  \ \ \ \text{and} \ \ \  \phi_c^2(x)=\sum_{n=-\infty}^{\infty}b_{n}e^{\frac{i\pi nx}{L}}.
\end{equation} 
Replacing the last equations in (\ref{ecuaordina}) we obtain that 
\begin{equation} \label{coefidelasuma}
ca_{n}[1+{\frac{\pi}{L}}n]-a_{n}=\frac{1}{2}\sum_{m=-\infty}^{\infty}a_{m}a_{n-m}, \ \  \forall n\in \mathbb{Z}.
\end{equation}
Now, from (\ref{ecuapoisson}) we consider $a_n=\frac{2\pi c}{L} e^{-\eta |n|},\ n\in\mathbb{Z}$ and $\eta > 0$. Substituting $a_n$ into the last identity we obtain 
\[\sum_{m=-\infty}^{\infty}a_{m}a_{n-m}=\tfrac{4\pi^2c^2}{L^2}e^{-\eta|n|}\left[|n|+1+2\sum_{k=1}^{\infty}e^{-2\eta k}\right] =\tfrac{4\pi^2c^2}{L^2}e^{-\eta|n|}(|n|+\coth\eta).\]
Therefore, we conclude that 
\begin{equation} \label{ultigual}
c\left[1 +{\frac{\pi}{L}}|n|\right]-1 = \frac{\pi c}{L}(|n|+\coth\eta),  \ \  \forall n\in \mathbb{Z}.
\end{equation}
We denote $\eta = \frac{w}{2(w-1)L}$ and consider $c \neq 1$ such that $0 < \frac{c}{c-1}<\frac{L}{\pi}.$ We choose $w=w(c)>1$ such that $\tanh (\eta) = \frac{\pi c}{(c-1)L}$, then we get from (\ref{ultigual}) that $\psi_{w(c)}=\phi_c$, hence, $\phi_c$ is given by (\ref{ecuapoisson}). Therefore, we obtain that $\phi_c$ has the form (\ref{solperiodicas}) with $\eta > 0$ satisfying  $\tanh (\eta) = \frac{\pi c}{(c-1)L}$.
\begin{remark}
$(1)$ Note from the fact that $c\neq 1$ satisfies $0 < \frac{c}{c-1}<\frac{L}{\pi}$ that we have three cases:\\
$(a)$ If $L=\pi$, then $c \in (-\infty,0)$\\
$(b)$ If $L<\pi$, then $c \in (1+\frac{\pi}{L-\pi},0)$\\
$(c)$ If $L>\pi$, then $c \in (-\infty,0) \ \cup \ (1+\frac{\pi}{L-\pi},+\infty)$\\
$(2)$ Observe that the sign of the solution $\phi_c$ depends on the sign of $c$, since we are interested in positive solutions (to apply the theory in \cite{anguloNatali}) we will suppose that $L>\pi$ and $c> 1 + \frac{\pi}{L-\pi}$.
\end{remark}

Also it is worth to note that if we consider $c=1$ in (\ref{ecuaordina}), then the unique real smooth solution that we obtain is $\phi\equiv 0.$ In fact, in this case we have that $\phi$ satisfies $\mathcal{H}\phi'-\frac{1}{2}\phi^2=0.$ Taking Fourier transform we arrive at
\[2|n|\widehat{\phi}(n)-\sum_{k=-\infty}^{+\infty}\widehat{\phi}(n-k)\widehat{\phi}(k)=0, \ \ \ \forall n\in\mathbb{N}.\]
In particular for $n=0,$ we have that $\sum_{k=-\infty}^{+\infty}\widehat{\phi}(-k)\widehat{\phi}(k)=0.$ Since $\phi$ is a real solution we have that $\widehat{\phi}(-k)=\overline{\widehat{\phi}(k)}.$ Therefore $\|\phi\|_{L^2_{per}}=0,$ then using the smoothness of $\phi$ we get that $\phi\equiv 0.$ \\

Now, using the fact that $\eta(c)=\tanh^{-1}\left(\frac{c\pi}{(c-1)L}\right)$ is a differentiable function if $c\neq 1$, we have the next result
\begin{prop}\label{procurva}
Let $L>\pi.$ Then the curve $c\in(1+{\frac{\pi}{L-\pi}} ,+\infty) \longrightarrow \phi_c \in H_{per}^{\frac 1{2}}([-L,L])$ is of class $C^1,$ where $\phi_c$ is given by (\ref{solperiodicas}). Furthermore, since $c>1+{\frac{\pi}{L-\pi}},$ we have that $\phi_c>0.$
\end{prop}
 

\section{Spectral Analysis}
This section is dedicated to study specific spectral properties associated to the linear operator $\mathcal{L} = c\mathcal{H} \partial_x -1 +c -\phi_c$, where $\phi_c$ is the periodic solution (\ref{solperiodicas}) given by proposition \ref{procurva} with fundamental period $2L,$ $L>\pi$ and $c>1+\frac{\pi}{L-\pi}$. This information will be basic in our stability theory for the rBO equation.\\

Our analysis will be on the periodic eigenvalues problem considered on $[-L,L]$\\ 
\begin{equation}\label{probautovalores}
\left \{ 
\begin{aligned}
\mathcal{L}\chi &= \lambda\chi  \\
\chi(-L) &= \chi(L), \ \  D^{\frac 1{2}}\chi(-L)=D^{\frac 1{2}}\chi(L). \\
\end{aligned} \right.
\end{equation}
We will show that problem (\ref{probautovalores}) determines exactly the existence of a single negative eigenvalue, which is simple, that zero is also an eigenvalue simple with eigenfunction $\phi_{c}'$ and the remainder of the spectrum is bounded away from zero. In this point we will used the theory developed by Angulo and Natali in \cite{anguloNatali}. In fact, initially deriving the equation (\ref{solitaria}) with regard to $x$ we get that $\mathcal{L} \phi_{c}' = 0$, therefore zero is an eigenvalue with eigenfunction associated $\phi_{c}' $. Additionally from the theory of compact self-adjoint operators we have that (\ref{probautovalores}) determines that the spectrum of $\mathcal{L}$ is a countable infinity set of eigenvalues $\{\lambda_n\}_{n=0}^{\infty}$ with
\[ \lambda_0 \leq \lambda_1 \leq \lambda_2 \leq ... \ ,\]
where $\lambda_n \rightarrow \infty$ as  $n \rightarrow \infty$ ( see for instance Proposition $3.1$ in \cite{anguloNatali}).\\

For the sake of completeness we will make here a summary of the theory given in \cite{anguloNatali}. In this work was studied the existence and the nonlinear stability of periodic travelling wave solutions for the family of equations\\
\begin{equation} \label{equaNatali1}
u_t+u^pu_x-(Mu)_x=0,
\end{equation}
where $p\geq 1$ is an integer and $M$ is a pseudo-differential operator in the context of periodic functions, defined through Fourier multipliers as 
\[\widehat{Mg}(k)=\zeta(k)\widehat{g}(k), \ \ \forall \ \  k\in \mathbb{Z},\]
where the symbol $\zeta$ of $M$ is a real, measurable, locally bounded and even function satisfying
\begin{equation} \label{desismb}
A_1|n|^{m_1}\leq \zeta(n)\leq A_2(1+|n|)^{m_2},
\end{equation} 
for $1\leq m_1\leq m_2,\  |n|\geq n_0, \ \zeta(n)>b$ for all $n\in \mathbb{Z}$ and $A_i>0, \ i=1,2.$
The main result in \cite{anguloNatali} is the determination of the spectrum of the linear, closed, non-bounded and self-adjoint operator $\mathcal{L}_0:D(\mathcal{L}_0) \longrightarrow L^2_{per}([-L,L])$ given by
\begin{equation} \label{operNatali}
\mathcal{L}_0u=(M+c)u-\phi^pu,
\end{equation}
where $D(\mathcal{L}_0)$ is dense in $L^2_{per}([-L,L])$ and $\varphi_c$ is a periodic travelling wave solution for the equation (\ref{equaNatali1}). The principal result reads as follows (see \cite{anguloNatali}).
\begin{theo}
Suppose that  $\varphi_c$ is a positive even solution of (\ref{equaNatali1}) such that  $\widehat{\varphi}_c > 0$ and $\widehat{\varphi_c^p} \in PF(2)$ discrete. Then\\
$(a)$ $\mathcal{L}_0$ has a unique negative eigenvalues $\lambda$, and it is simple;\\
$(b)$ the eigenvalue $0$ is simple.\\
\end{theo}
Here a sequence $\alpha = (\alpha_n)_{n \in \mathbb{Z}} \subseteq \mathbb{R}$ is in the class $PF(2)$ discrete if
\begin{equation}\label{pf2} 
\begin{aligned}
&(i)\ \alpha_n > 0 \ \  \text{for all} \ \  n \in \mathbb{Z}, \\
&(ii)\ \alpha_{n_1-m_1}\alpha_{n_2-m_2}-\alpha_{n_1-m_2}\alpha_{n_2-m_1}>0 \  \text{for} \ \  n_1 < n_2 \ \ \text{and} \ \  m_1 < m_2.
\end{aligned}
\end{equation}
Although the rBO does not have the form (\ref{equaNatali1}), we still can apply their result to get the required information about the spectrum of the operator in (\ref{operadorMio}). For this, consider $L>\pi,$   $c>1+\frac{\pi}{L-\pi}$ and define $ M = c \mathcal{H} \partial x  -1,$ then $\mathcal{L} \ = (M + c) - \phi_{c}.$ In this case,
\[ \widehat{Mf}(k) = \zeta (k) \widehat{f}(k), \ \  \forall  k  \in  \mathbb{Z}\] 
where $\zeta(k) \equiv c|k| - 1.$ Furthermore, considering $A_{1}=1,$  $A_2=c$ and $m_1=m_2=1,$  there exists  $N_{0} \in \mathbb{N}$ such that $\zeta$ satifies (\ref{desismb})  for all $k\geq N_0.$\\

Since $c>1+{\frac{\pi}{L-\pi}},$ and $\eta > 0$, then  $\phi_{c} > 0$ and  it is clear that  $a_n=\widehat {\phi_{c}}(n)= \frac{4\pi^2}{L^2} e^{-\eta|n|} > 0$.\\

Before to continue with the study of the spectrum of $\mathcal{L}$ given in (\ref{operadorMio}), we note that the main theorem in \cite{anguloNatali} is still valid if we replace $(ii)$ in (\ref{pf2}) by the weaker condition\\
\[ 
(ii') \left \{ 
\begin{aligned}
\alpha_{n_1-m_1}\alpha_{n_2-m_2}-\alpha_{n_1-m_2}\alpha_{n_2-m_1}\geq 0&,\  \text{for all} \  n_1<n_2 \ \text{and} \ m_1<m_2  ,\\
\alpha_{n_1-m_1}\alpha_{n_2-m_2}-\alpha_{n_1-m_2}\alpha_{n_2-m_1}>0&,\   \text{if} \  n_1<n_2 ,\  m_1<m_2, \ n_2>m_1,\\& \ \ \text{and}\ n_1<m_2.  \\
\end{aligned} \right.\]
So, in order to use the result obtained in \cite{anguloNatali}  we just have to show that $a_n=\frac{2c\pi}{L} e^{-\eta|n|}$ satisfies $(ii'),$ which is equivalent to prove that   
\[(a)\ |n_1-m_1|+|n_2-m_2| \leq |n_1-m_2|+|n_2-m_1|,\  \text{if} \  n_1<n_2 \ \text{and} \ m_1<m_2, \ \text{and}\] 
\[(b)\ |n_1-m_1|+|n_2-m_2| < |n_1-m_2|+|n_2-m_1|,\   \text{if} \  n_1<n_2 ,\  m_1<m_2, \ n_2>m_1\] \[\hspace{5.5cm}\text{and}\ n_1<m_2. \]
The proof of $(a)$ and $(b)$ are easy, so we skip the details, then  as a consequence of this analysis we have the next result.
\begin{prop}\label{prospec}
Let $\phi_c$ be the periodic wave solution given by Proposition \ref{procurva}, with $c>1+\frac{\pi}{L-\pi}$ and 
$L>\pi.$ Then, the linear operator $\mathcal{L}$ define by (\ref{operadorMio}) with domain   $H_{per}^\frac{1}{2}([-L,L])\subseteq L^2_{per}([-L,L])$ has its first two eigenvalues simple with zero being the second one (with eigenfunction $\frac{d}{dx}\phi(x)$). Moreover, the remainder of the spectrum is constituted by a discrete set of eigenvalues which converge to $+\infty.$ 
\end{prop}


\section{Stability of travelling wave solutions}
In this section we will use the Lyapunov method for studying the nonlinear stability of solutions $u(x,t)=\phi_c(x-ct)$ with $\phi_c$ given by Proposition \ref{procurva}. For this, the use of the conservations laws for (\ref{ecuabasica})
\[E(u)= \frac{1}{2}\int_{-L}^{L} ((D^{\frac 1{2}}u)^2 -\frac{1}{3} u^3)dx\ \ \text{and}\ \  F(u)\equiv \frac{1}{2}\int_{-L}^{L} (u^2+(D^{\frac 1{2}}u)^2)dx\] 
will be required. The notion of stability that we will prove is the \textit{orbital stability}, more precisely, we shall prove that the orbit generated by $\phi_c$,
\[\mathcal{O}_{\phi_c}=\left\{\phi_c(\cdot+y): y\in \mathbb{R} \right\}\]
is stable for the periodic flow generated by the rBO equation. Namely, for every $\epsilon>0$ there exists $\delta(\epsilon)>0$ such that if $u_0\in H^{\frac 3{2}}_{per}$ and
\[\inf_{y\in \mathbb{R}}\|u_0-\phi_c(\cdot+y)\|_{H^{\frac 1{2}}_{per}}<\delta,\]
then the solution $u$ of the rBO equation (\ref{ecuabasica}) with initial datum $u_0$ satisfies 
\[\|u(t)-\phi_c(\cdot+y)\|_{H^{\frac 1{2}}_{per}}<\epsilon,\]
for all $t\in \mathbb{R}$ and $y=y(t).$\\
Note that we do not have a well-posedness result in $H_{per}^{\frac 1{2}},$ for this reason in our definition  of stability we took initial data $u_0$ in $H^{\frac 3{2}}_{per}.$ \\

Before to establish our main result we will prove an useful Lemma.

\begin{lemm} \label{lema4}
Let $\phi_c$ be the wave solution given by Proposition \ref{procurva} with $c\in\left(1+\frac{\pi}{L-\pi},+\infty\right)$ and $L>\pi.$ Then, the linear operator $\mathcal{L}= c \mathcal{H} \partial x -1+c-\phi_{c}$ satisfies that
\begin{align}
(a)&\ \  \alpha:=\inf \{(\mathcal{L}f,f) :\|f\|_{L_{per}^2} = 1 \ \text{and}\ \ (f,\phi_c+\mathcal{H}\phi_c')=0\}= 0,\label {igualalpha} \\
(b)&\ \  \beta:=\inf \{(\mathcal{L}f,f):\|f\|_{L_{per}^2}=1, (f,\phi_c+\mathcal{H}\phi_c')=0 \ \  \text{and} \ \ (f,\phi_c\phi_c')=0\}>0 \label{igualpos}
\end{align}
\end{lemm}
\proof
(a) Because $\phi_c$ is bounded, it is inferred that $\alpha$ is finite. Since $(\phi,\phi_c+\mathcal{H}\phi_c')=0$ and $\mathcal{L}\phi_c' = 0$ it follows that $ \alpha \leq 0.$ Next we will show that the inf in (\ref{igualalpha}) is attained. In fact, since $\alpha$ is finite, there exists a sequence
$\{f_j\}^\infty_{j=0} \subset H_{per}^{\frac{1}{2}} $ with $\|f_j\|_{L_{per}^2} = 1,( f_j,\phi_c+\mathcal{H}\phi_c' ) = 0 $  and $(\mathcal{L}f_j,f_j)\rightarrow \alpha$ as $j\rightarrow \infty.$  It follows that $\|f_j\|_{H_{per}^\frac{1}{2}}$ is uniformly bounded as $j$ varies. So, there exists a subsequence of $f_j,$ which we denote $\{f_j\}$ again and a function $ f^* \in  H_{per}^{\frac{1}{2}}$ such that $ f_j \rightharpoonup f^* $ in $ H_{per}^{\frac{1}{2}}.$ Now, since the embedding $H^{\frac 1{2}_{per}} \hookrightarrow L^2_{per}$ is compact we obtain that  $(f^*,\phi_c+\mathcal{H}\phi_c')=0$  and
$(\phi_c f_j,f_j) \longrightarrow  (\phi_c f^*,f^*)$ when $ j\to\infty.$ So $f^* \neq 0$ and since the weak convergence is lower continuous we obtain
\[\|D_x^{\frac{1}{2}}f^*\|^2_{L_{per}^2} \leq \ \liminf_{j\to\infty}\|D_x^{\frac{1}{2}}f_j\|^2_{L_{per}^2}.\]  
Now, define $f=\frac{f^*}{\|f^*\|_{L^2} },$ then $(f,\phi_c+\mathcal{H}\phi_c')=0,$   $\|f\|_{L_{per}^2}=1$ and
\[\alpha\leq(\mathcal{L}f,f)\leq\frac{\alpha}{\|f\|_{L_{per}^2}^2}\leq\alpha.\]
Therefore, $\alpha$ is a minimum.
Now, we want to show that $\alpha \geq 0$. In this case, we will apply the Lemma E1 in \cite{weinstein2} ( which works in the periodic setting ) in the case that $A =\mathcal{L}$ and $R=\phi_c+\mathcal{H}\phi_c'.$  In fact, from Proposition \ref{prospec}, $\mathcal{L}$ has the spectral properties required by Lemma E1. Next, we need to find  $\chi$ such that $\mathcal{L}\chi = \phi_c+\mathcal{H}\phi_c'$ and $(\chi,\phi_c+\mathcal{H}\phi_c') \leq 0$. In fact, from Proposition \ref{procurva} we have that the mapping  $c\in(1+{\frac{\pi}{L-\pi}} ,+\infty) \longrightarrow \phi_c \in H_{per}^{\frac 1{2}}([-L,L])$ is of class $C^1,$ so by differentiating (\ref{ecuaordina}) with regard to $c$ we obtain that $\chi= -\frac{d}{dc}\phi_c$ satisfies 
\[\mathcal{L}(\chi)=\phi_c+\mathcal{H}\phi_c'.\]
Observe that,
\begin{align*}
(\chi,\phi_c+\mathcal{H}\phi_c')&=-L\frac d{dc}\sum_{n=-\infty}^{\infty}(1+|n|)|\widehat{\phi_c}(n)|^2=-\frac{8\pi^2c}{L}
\sum_{n=-\infty}^{\infty}(1+|n|)e^{-2|n|\eta}\\&+\frac{8\pi^2c^2}{L}\frac{d\eta}{dc}
\sum_{n=-\infty}^{\infty}(1+|n|)|n|e^{-2|n|\eta}.
\end{align*}
But, from (\ref{valortanh})  and using the fact that $c>1+\frac{\pi}{L-\pi},$ we have that 
\begin{align*}
\frac{d\eta}{dc}&=\frac{d}{dc} \left( \tanh^{-1}\left(\frac{c\pi}{(c-1)L}\right)\right)=-\frac{\pi}{(c-1)^2L}\left(1-\left( \frac{c\pi}{(c-1)L}\right)^2\right)^{-1}<0
\end{align*}
Therefore $(\chi,\phi_c+\mathcal{H}\phi_c') < 0$ and  Lemma E1 give us that $\alpha \geq 0.$ This finish the proof of (\ref{igualalpha}).\\

For part (b).  From part (a) is inferred that $\beta \geq 0$. Suppose that $\beta = 0.$ then we can find a function $f$ such that $\|f\|_{L_{per}^2} = 1$ and $(f,\phi_c+\mathcal{H}\phi_c')=(f,\phi_c\phi_c')=(\mathcal{L}f,f)= 0.$ Therefore, there exist  $\gamma,\theta,\nu$ such that
\[\mathcal{L}f = \gamma f + \theta(\phi_c+\mathcal{H}\phi_c') + \nu\phi_c\phi_c'.\]
So, $\gamma=\nu=0.$ Therefore $\mathcal{L}f = \theta(\phi_c+\mathcal{H}\phi_c').$ Now consider  $\chi= -\frac{d}{dc}\phi_c,$ it follows that $\mathcal{L}(f-\theta\chi)=0,$ then $(f-\theta\chi,\phi_c+\mathcal{H}\phi_c')=0=(f,\phi_c+\mathcal{H}\phi_c')-\theta(\chi,\phi_c+\mathcal{H}\phi_c'),$  so, $\theta = 0,$ because $(\chi,\mathcal{H}\phi_c')\neq 0,$ therefore $\mathcal{L}f=0,$ and so there exists a  $\lambda \in \mathbb{R}-\{0\}$ such that  $f=\lambda\phi_c',$ and hence $f$ is orthogonal to $\phi_c\phi_c',$ which is a contradiction. Therefore $\beta > 0$ and the proof of the lemma is completed.
\endproof

We note that from (\ref{igualpos}) and from the specific form of $\mathcal{L}$ we have that if  $(f,\phi_c+\mathcal{H}\phi_c')= 0$  and $(f,\phi_c\phi_c')=0,$ then there exists $\beta_{0}>0$ such that
\[(\mathcal{L}f,f) \geq \beta_{0}\|f\|^2_{H_{per}^\frac{1}{2}}.\]
Now we establish our main result.


\begin{theo} \label{teoremaPrinc1}
Let  $L>\pi$  and $\phi_c$ be the periodic wave solution given by Proposition \ref{procurva} with $c\in \left(1+\frac{\pi}{L-\pi},+\infty\right).$  Then the orbit $\mathcal{O}_{\phi_c}$ is nonlinear stable with regard to the periodic flow generated by the rBO equation.
\end{theo}
\proof
The proof is based in the ideas developed in \cite{benjamin1}, \cite{bona2}, \cite{weinstein3}, \cite{anguloNatali}. We shall give only an outline of the proof. Initially, note that  $F'(u)=u +\mathcal{H}u_x$ and $E'(u)=\mathcal{H}u_x-\frac{1}{2}u^2$, then $\phi_c$ is a critical point of the functional $\mathcal{B} \equiv E + (c-1)F$. Additionally, since $F''(u)=1 +\mathcal{H}\partial_x$ and $E''(u)=\mathcal{H}\partial_x-u$, we have 
\[ E''(\phi_c)+(c-1)F''(\phi_c)=c\mathcal{H}\partial_x +(c-1) -\phi_c= \mathcal{L}.\]
Now, define for $r\in[-L,L]$ and $t\in \mathbb{R},$
\[\Omega_t(r)\equiv\|D^\frac{1}{2}u(\cdot +r,t)-D^\frac{1}{2}\phi_c\|_{L_{per}^2}^2+\tfrac{c-1}{c}\|u(\cdot +r,t)-\phi_c\|_{L_{per}^2}^2.\]  
then, using standard arguments (see \cite{benjamin1,bona2}) there exists an interval of time $I=[0,T]$ such that the $\inf_{r\in\mathbb{R}}\Omega_t(r)$ is attained in $\gamma=\gamma(t)$ for every $t\in I.$ Hence, we get that
\begin{equation}\label{inf}
\Omega_t(\gamma(t))=\inf_{r\in\mathbb{R}}\Omega_t(r).
\end{equation}
Now, consider the perturbation of the periodic travelling wave $\phi_c$
\begin{equation}\label{pertur}
u(x+\gamma ,t)=\phi_c (x)+v(x,t)
\end{equation} 
for $t\in[0,T]$ and $\gamma=\gamma(t)$ determined by (\ref{inf}). Then, differentiating $\Omega_t(r)$ with respect to $r,$ evaluating at values that minimize $\Omega_t(r)$ and using (\ref{pertur}) we obtain that $v$ satisfies the compatibility relation 
\begin{equation} \label{condigual0}
\int_{-L}^{L}\phi_c'(x) \phi_c(x) v(x,t)dx=0,
\end{equation}
for all $t\in[0,T].$ 
Next, using the fact that $E$ and $F$ are conserved quantities, the representation (\ref{pertur}), the embedding $H_{per}^{\frac 1{2}}([-L,L])\hookrightarrow L^r([-L,L]),$ for all $r\geq 2$ and the fact that $\phi_c$ satisfies (\ref{ecuaordina}), we conclude 
\begin{align}\label{suitbound}
\Delta \mathcal{B}(t)&=\mathcal{B}(u_0)-\mathcal{B}(\phi_c)=\mathcal{B}(u(\cdot,t))-\mathcal{B}(\phi_c)=\mathcal{B}(\phi_c +v(\cdot,t))-\mathcal{B}(\phi_c)\notag \\& \geq \frac{1}{2}(\mathcal{L}v,v) -c_0\|v\|_{H_{per}^\frac{1}{2}}^3.
\end{align}
where $c_0$ is a positive constant. To obtain our result we need to establish a suitable bound for the quadratic form in (\ref{suitbound}). Initially, we consider the normalization  $F(u_0)=F(\phi_c),$ then, \[\int_{-L}^{L} u^2(t)+(D^{\frac 1{2}}u(t))^2 \ dx = \int_{-L}^{L} \phi_c^2+(D^{\frac 1{2}}\phi_c)^2 \ dx\] 
for all  $t \in [0,T].$ By (\ref{pertur}) it follows 
\[-2(v,\phi_c+\mathcal{H}\phi_c')=\|v(t)\|_{L_{per}^2} +\|D^{\frac{1}{2}}v\|_{L_{per}^2}.\]
Without loss of generality, we suppose that  $\|\phi_c+\mathcal{H}\phi_c'\|_{L_{per}^2}=1.$ Define $v_{\parallel}$ and  $v_{\perp}$ as $v_{\parallel}=(v,\phi_c+\mathcal{H}\phi_c') (\phi_c+\mathcal{H}\phi_c')$ and $v_{\perp}= v-v_{\parallel}.$ So, $(v_{\perp},\phi_c+\mathcal{H}\phi_c')=0$  and $(v_{\perp},\phi_c\phi_c')=0.$
By (\ref{condigual0}) and Lemma \ref{lema4} it follows that
\begin{align} \label{desiformcuad1}
(\mathcal{L} v_{\perp},v_{\perp}) &\geq \beta \|v_{\perp}\|^2_{L_{per}^2}\geq\beta \|v\|_{L_{per}^2}^2-\tilde{\beta_3}\|v\|^4_{H^{\frac 1{2}}_{per}}
\end{align}
with $\beta,\beta_3 >0.$ Again, without loss of generality suppose that $(\mathcal{L}(\phi_c+\mathcal{H}\phi_c'),\phi_c+\mathcal{H}\phi_c') <0,$ then 
\begin{align} \label{desiformcuad2}
(\mathcal{L}v_{\parallel},v_{\parallel})\geq-\tilde{\beta_4}\|v\|^4_{H^{\frac 1{2}}_{per}}.
\end{align}
Furthermore, using the Cauchy-Schwarz inequality we get
\begin{align} \label{desiformcuad3}
(\mathcal{L}v_{\parallel},v_{\perp})&\geq-\tilde{\beta_{2}}\|v(t)\|^3_{H_{per}^\frac{1}{2}},
\end{align}
where $\tilde{\beta_{j}} > 0,$ for $  j=3,4.$\\
Now, using (\ref{desiformcuad1}), (\ref{desiformcuad2}), (\ref{desiformcuad3}) and the specific form for $\mathcal{L}$ we conclude that
\begin{equation} \label{desimportatisima}
(\mathcal{L}v,v)\geq\beta_{0} \|v(t)\|^2_{H_{per}^\frac{1}{2}}-\beta_1\|v(t)\|^3_{H_{per}^\frac{1}{2}}-\beta_{2} \|v(t)\|^4_{H_{per}^\frac{1}{2}},
\end{equation}
where $\beta_j>0,$ for $j=0,1,2.$  Hence, from (\ref{inf}), (\ref{suitbound}) and  (\ref{desimportatisima}) it follows that for all $t\in [0,T]$ 
\begin{align}\label{desiend}
\Delta \mathcal{B}(t)\ge g\left(\|v(t)\|_{\frac{1}{2},c}\right),
\end{align} 
where   $\|f\|^2_{\frac{1}{2},c}:=\|D^\frac{1}{2}f\|^2_{L^2_{per}}+ \tfrac{c-1}{c}\|f\|^2_{L^2_{per}}$ and  $g(s)=s^2 -\sum_{k=3}^{4}d_k(c)s^k,$  with  $\eta, d_k>0.$ The essential properties of $g$ are $g(0)=0$ and  $g(s)>0$ for $s$ small. The stability result is an immediately consequence of (\ref{desiend}). In fact, let  $\epsilon>0$ small enough such that $g(\epsilon)>0.$ Then by using the properties that $\mathcal{B}$ is uniformly continuous on $S:= \left\{u\in H_{per}^{\frac{1}{2}}: F(u)=F(\phi_c)  \right\},$ $\Delta \mathcal{B}(t)$ is constant in time and $t\mapsto\|v(t)\|^2_{\frac{1}{2},c}$ is a continuous function, we have that there is  $\delta(\epsilon)>0$ such that if $v\in S$  and $\|v-\phi_c\|_{\frac{1}{2},c}< \delta$ then for $t\in[0,T],$ 
\begin{equation}\label{finalineq}
g\left(\|v(t)\|_{\frac{1}{2},c}\right)\leq \Delta \mathcal{B}(0)\leq|\Delta \mathcal{B}(0)|<g(\epsilon)  \Rightarrow \|v(t)\|^2_{\frac{1}{2},c}<\epsilon.
\end{equation} 
Which shows that $\mathcal{O}_{\phi_c}$ is orbitally stable in $H^{\frac 1{2}}_{per}([-L,L])$ relative  to small perturbations which preserve the $L_{per}^2$ norm.\\
The inequality (\ref{finalineq}) is still true for $t>0,$ this is an immediately consequence of the fact that the the mapping $t \mapsto\inf_{r\in\mathbb{R}}\Omega_t(r)$ is continuous (see \cite{bona2}).\\

To prove stability to general perturbations, we use that the mapping $c\in(1+{\frac{\pi}{L-\pi}} ,+\infty) \mapsto\phi_c \in H_{per}^{\frac 1{2}}([-L,L])$ is continuous, that the mapping $c\in(1+{\frac{\pi}{L-\pi}} ,+\infty)\mapsto F(\phi_c)$ is strictly increasing, the preceding theory and the triangle inequality (see \cite{bona2}, \cite{weinstein3}, \cite{angulo1}). Then Theorem \ref{teoremaPrinc1} is proved.
\endproof

\section{Generalization of the theory}
In this section we will extend the theory developed for the rBO to a family of equations of the form (\ref{formagen}). $H$ is defined as a Fourier multiplier by
\[\widehat{Hu}(n)=\alpha(n)\widehat{u}(n), \ \ \ n\in\mathbb{Z},\]
where the symbol $\alpha$ is assumed to be a real, mensurable, locally bounded, even function on $\mathbb{R}$, satisfying the conditions
\begin{equation}\label{desSymGen}
A_1|n|^{m_1}\leq\alpha(n)\leq A_2(1+|n|)^{m_2},
\end{equation}
where $1\leq m_1\leq m_2,\ |n|>n_0, \ \alpha(n)>b$ for all $n\in\mathbb{Z}$  and $A_i>0,$ for $i=1,2.$\\

The travelling wave solutions in which we are interested will have again the form
\begin{equation}\label{formaSolGen}
u(x,t)=\phi(x-ct),
\end{equation}
where the profile $\phi:\mathbb{R}\rightarrow \mathbb{R}$ is a smooth periodic function with fundamental period (\textit{a priori}) $2L>0.$\\

Then substituting the form (\ref{formaSolGen}) in (\ref{formagen}) and integrating once (and considering the constant of integration zero in our theory) we arrived at
\begin{equation}\label{EDOGen}
cH\phi_c+(c-1)\phi_c-\frac 1{p+1}\phi_c^{p+1}=0.
\end{equation}
As it is well known the equation (\ref{formagen}) has two conservation laws
\[E(u)=\frac 1{2}\int_{-L}^{L}uHu-\frac{2}{(p+1)(p+2)}u^{p+2}\ dx\]
and 
\[F(u)=\frac 1{2}\int_{-L}^{L}uHu+u^2\ dx.\]
If we consider \textit{a priori} the existence of a periodic travelling wave solution $\phi_c$ for the equation (\ref{EDOGen}), it is easy to see that $E'(\phi_c)+(c-1)F'(\phi_c)=0.$\\

Now, define $\mathcal{L}:=E''(\phi_c)+(c-1)F''(\phi_c)=cH+(c-1)-\phi_c^p,$ that is,
\begin{equation}\label{operGen}
\mathcal{L}u=cHu+(c-1)u-\phi_c^pu.
\end{equation}
Then the  operator $\mathcal{L}:D(\mathcal{L})\rightarrow L^2_{per}([-L,L])$ is  linear, closed, not bounded and self-adjoint defined on a dense subset of $ L^2_{per}([-L,L]).$ Also it is easy to see that $\mathcal{L}\phi_c'=0,$ this is, zero is a eigenvalue of $\mathcal{L}$ with eigenfunction $\phi_c'.$\\

Following the proof that we already made to get the orbital stability of periodic traveling wave solutions  for the rBO equation, the conditions which prove stability are
\begin{equation}\label{condGen}
\begin{array}{cl}
(C_0)& \text{there is a nontrivial smooth curve of periodic solutions for (\ref{EDOGen})}\\
    & \text{of the form}\  c\in I\subset\mathbb{R}\rightarrow \phi_c\in H^{m_2}_{per}([-L,L]);\\
(C_1)& \mathcal{L}\  \text{has an unique negative eigenvalue and it is simple};\\
(C_2)& \text{the eigenvalue zero is simple};\\
(C_3)& \frac{d}{dc}\int_{-L}^{L}\phi_cH\phi_c+\phi_c^2\ dx>0.
\end{array}
\end{equation}

Next we will show that is possible to use the theory established in \cite{anguloNatali} to give suficient conditions to obtain $C_1$ and $C_2$ for the operator $\mathcal{L}$ associated to the problem (\ref{formagen}). Let us start with a proposition which says that the essential spectrum of $\mathcal{L}$ is empty. Its proof follows the same ideas as in Proposition 3.1.1 in \cite{anguloNatali}.
\begin{prop}
The operator $\mathcal{L}$ in (\ref{operGen}) is closed, unbounded, self-adjoint on $L^2_{per}([-L,L])$
whose spectrum consists of an enumerable (infinite) set of eigenvalues which converge to $+\infty$. In particular, $\mathcal{L}$ has
zero as an eigenvalue with eigenfunction $\frac{d}{dx}\phi_c.$
\end{prop}

Next we establish the principal result of this section.

\begin{theo}\label{teorPrinGen}
Let $\phi_c$ be a positive even solution of (\ref{EDOGen}). Assume that $\widehat{\phi_c} > 0$ and $\widehat{\phi_c^p}\in PF(2)$ discrete, then $(C_1)$ and $(C_2)$ in (\ref{condGen}) hold for the
operator $\mathcal{L}$ in (\ref{operGen}).
\end{theo}
\proof
 Note that the operator $\mathcal{L}$ can be written as 
\[\mathcal{L}u=(M+c)u-\phi_c^pu,\]
where $M=cH-1.$ The symbol of $M$ is $\zeta(n)=c\alpha(n)-1.$ Using (\ref{desSymGen}) it is easy to see that for all $c\neq 0$ there is $n_0\in\mathbb{N}$ such that
\[B_1|n|^{m_1}\leq\zeta(n)\leq B_2(1+|n|)^{m_2},\ \ \ \ \forall n\geq n_0,\]
where $B_1=\frac{2}{|c|A_1}$ and $B_2=cA_2+1.$ Then we can apply Theorem 4.1 in \cite{anguloNatali} to obtain that $C_1$ and $C_2$ hold for the operator $\mathcal{L}.$
\endproof
\subsection{Stability of the BBM equation}
Next, we are interested in applying the Theorem \ref{teorPrinGen} to obtain the orbital stability of travelling wave solutions associated to the BBM equation. Let us start with the definition of stability in the case of equation (\ref{formagen}).
\begin{defi}
Let $\phi_c$ be a periodic travelling-wave solution with period $2L$ of (\ref{EDOGen}). We define the set 
$\Omega_{\phi}\subset H^{\frac{m_2}{2}}_{per},$ the orbit generated by $\phi$, as
\[\Omega_{\phi}=\{f:f=\phi(\cdot+r) \  \text{for some} \ r\in \mathbb{R}\}.\]
And, for any $\gamma> 0$, we define the set $U_{\gamma}\subset H^{\frac{m_2}{2}}_{per}$ by
\[U_{\gamma}=\{f: \inf_{g\in\Omega_{\phi}}\|f-g\|_{ H^{\frac{m_2}{2}}_{per}}<\gamma\}.\]
With this terminology, we say that $\phi_c$ is (orbitally) stable in $H^{\frac{m_2}{2}}_{per}$ by the flow
generated by (\ref{formagen}) if the following hold:\\

$(i)$ There is $s_0$ such that $H^{s_0}_{per}\subset H^{\frac{m_2}{2}}_{per}$ and the initial value problem associated to (\ref{formagen}) is globally well-posed in $H^{s_0}_{per}$.\\

$(ii)$ For every $\epsilon> 0$, there is $\delta > 0$ such that, for all $u_0\in U_{\delta}\cap H^{s_0}_{per}$, the solution
$u$ of (\ref{formagen}) with $u(0, x) = u_0(x)$ satisfies $u(t)\in U_{\epsilon}$ for all $t > 0$.
Otherwise, we say that $\phi_c$ is unstable in $H^{\frac{m_2}{2}}_{per}$.
\end{defi}

The proof of the following general stability theorem can be shown by following the ideas used in the rBO stability theorem.
\begin{theo}\label{teoremI}
Let $\phi_c$ be a periodic travelling-wave solution of (\ref{EDOGen}), and suppose
that part $(i)$ of the definition of stability holds. Suppose also that the operator $\mathcal{L}$ defined previously in (\ref{operGen}) has properties $(C_1)$ and $(C_2)$ in (\ref{condGen}). Choose  $\chi \in  L_{per}^2$ such that $\mathcal{L}\chi = \phi_c+H\phi_c$, and define $I=(\chi,\phi_c+H\phi)_{L_{per}^2}$. If $I<0$ then $\phi_c$ is stable.
\end{theo}

\begin{remark}
In our cases the function $\chi$ in Theorem \ref{teoremI} is $\chi = -\frac{d}{dc}\phi_c$.
\end{remark}

Now, we apply the results obtained previously to the proof of the stability of
periodic travelling-wave solutions of cnoidal type associated with the BBM equation
\begin{equation}\label{equaBBM}
u_t+u_x+uu_x-u_{xxt}=0.
\end{equation}
In this case the travelling wave solutions associated to (\ref{equaBBM}) satisfy
\begin{equation}\label{edoBBM}
c\phi_c^{''}-(c-1)\phi_c + \frac{1}{2} \phi_c^2 = 0.
\end{equation}
We consider the solitary wave solutions
\[\varphi_w(x)=3(w-1)\text{sech}^2\left(\sqrt{\frac{w-1}{w}} \frac{x}{2}\right),\]
with $w>1$, whose Fourier transform is given by \[\widehat{\varphi}_w^{\mathbb{R}}(\xi)=12\pi\xi w\hspace{1pt} \text{csch}\left(\sqrt{\frac{w}{w-1}}\pi\xi\right).\]
Then from the Poisson summation theorem we consider
\begin{align*}
\psi_{w}(\xi) &=\frac{12w}{L} \sqrt{\tfrac{w-1}{w}} + \frac{24 \pi w}{L^2} \sum_{n=1}^{\infty} n \ \text{csch} \left(\sqrt{\tfrac{w}{w-1}}\frac{\pi n}{L}\right)\cos \left(\frac{2 \pi n \xi}{L} \right).
\end{align*}

Note that for each $L>0$ fixed, there exists a unique $k_0\in(0,1)$ such that  $\frac{K(k)}{K(k')L}<1$ and since $w>1$ is arbitrary, we can consider $w := w(k)$ such that $\sqrt{\frac{w-1}{w}}= \frac{K(k)}{K(k')L},$ for all $k\in(0,k_0),$ where ${k'}^2 = 1-k^2.$ Then, we obtain
\begin{equation}\label{solPoissonBBM}
\psi_{w(k)}(\xi)=\frac{12w}{L}\sqrt{\tfrac{w-1}{w}}+\frac{24\pi w}{L^2}\sum_{n=1}^{\infty}n \ \text{csch} \left(\frac{\pi n K'}{K}\right)\cos \left(\frac{2 \pi n \xi}{L} \right).
\end{equation}
Now, we use the Fourier expansion of $\text{dn}^2$ (see \cite{oberhettinger11}), that is,
\[K^2\left[\text{dn}^2 \left(\frac{2K\xi}{L};k\right)-\frac{E}{K}\right]= 2\pi\sum_{n=1}^{\infty}\frac{nq^n}{1-q^{2n}}\cos \left( \frac{2 \pi n \xi}{L} \right),\]
where $q=e^{-{\left(\frac{\pi K'}{K}\right)}}.$ We can conclude that
\[\frac{q^{n}}{1-q^{2n}}=\frac{1}{2}\ \text{csch} \left( \frac{n \pi K'}{K}\right).\]
Thus we obtain from (\ref{solPoissonBBM}) that
\begin{equation}
\psi_{w(k)}(\xi) = \frac{12w}{L}\sqrt{\tfrac{w-1}{w}} + \frac{24 K^2 w}{L^2} \left[\text{dn}^2 \left(\frac{2K\xi}{L};k\right)-\frac{E}{K}\right],
\end{equation}
for $k \in (0,k_0)$ fixed.\\

Next, because the last equality we will consider $\phi_c(x)= a+b\left[\text{dn}^2 (d \xi;k)-\frac{E}{K}\right]$ a
periodic travelling-wave solution for (\ref{edoBBM}) with period $L$. Then, the following nonlinear
system is obtained:
\[\left\{
\begin{array}{l}
 \frac{b^2}{2} - 6cbd^2 = 0\\
4bd^2c(1+{k'}^2)+ab -b^2\frac{E}{K} - (c-1)b  = 0\\
\frac{a^2}{2}-ab\frac{E}{K} +\frac{b^2}{2}(\frac{E}{K})^2-(c-1)a+(c-1)b\frac{E}{K}-2cbd^2{k'}^2 = 0.
\end{array}\right.\]
Since $\phi_c$ is periodic of period $L$ it follows that $d=\frac{2K}{L}.$ Then, from the first equation of the system above we have that $b=\frac{48cK^2}{L^2}.$ Substituting those values at the second equation we get
\[aL^2 = 48cKE - 16K^2c(2-k^2)+cL^2 -L^2 = 0.\]
Plugging the  values of $aL^2$, $d$ and $b$ in the third equation of the system we arrived at the following quadratic equation
\[\left[256K^4(1-k^2+k^4)-L^4\right]c^2+2cL^4-L^4 = 0.\]
The solutions of this equation are
\[c=\frac{L^2}{L^2+16K^2\sqrt{1-k^2+k^4}} \ \ \ \text{and} \ \ \ \  c=\frac{L^2}{L^2-16K^2\sqrt{1-k^2+k^4}}.\]

Next, if we want to use Theorem \ref{teorPrinGen}  we need that $cH-1+c>0$, that is, $ck^2-1+c>0$ for all $k \in \mathbb{Z}$, in particular for $k=0$ we obtain the $c>1$. If we choose $c$ as
\[c=\frac{L^2}{L^2+16K^2\sqrt{1-k^2+k^4}}\]
we will have that $0<c<1$, then we cannot apply Theorem \ref{teorPrinGen} in this case. For this reason, for $L>2\pi$ fixed, we will choose
\[c=\frac{L^2}{L^2-16K^2\sqrt{1-k^2+k^4}},\ \ \ \ \forall k\in(0,k_L),\]
where $k_L\in(0,1)$ is such the $L^2-16K^2\sqrt{1-k^2+k^4}>0$ for all $ k\in(0,k_L).$\\

Observe that for $L>2\pi$ fixed, we get that $k\mapsto c(k)$ is an increasing function on $(0,k_L)$ (see Figure 1), therefore for all $k \in (0,k_L)$ we have $c \in \left(c^{\ast}, +\infty \right),$ where $c^{\ast}=1+\frac{4\pi^2}{L^2-4\pi^2}.$\\

We can write $\phi_c$ in terms of $\text{cn}^2$, that is
\begin{equation}\label{curvaBBM}
\phi_c(x)=\beta_2 + (\beta_3-\beta_2) \ \text{cn}^2 \left( \sqrt{\frac{\beta_3-\beta_1}{12c}} \ x; k\right),
\end{equation}
where
\[\beta_2 =\frac{16cK^2(2{k'}^2-1)}{L^2} +c-1, \ \ \ \beta_3=\frac{16cK^2}{L^2}(1+k^2) + c-1\]
and $\beta_1$ is such that
\[\beta_3-\beta_1 =\frac{48cK^2}{L^2}.\]
\includegraphics[scale=0.29, angle=-90]{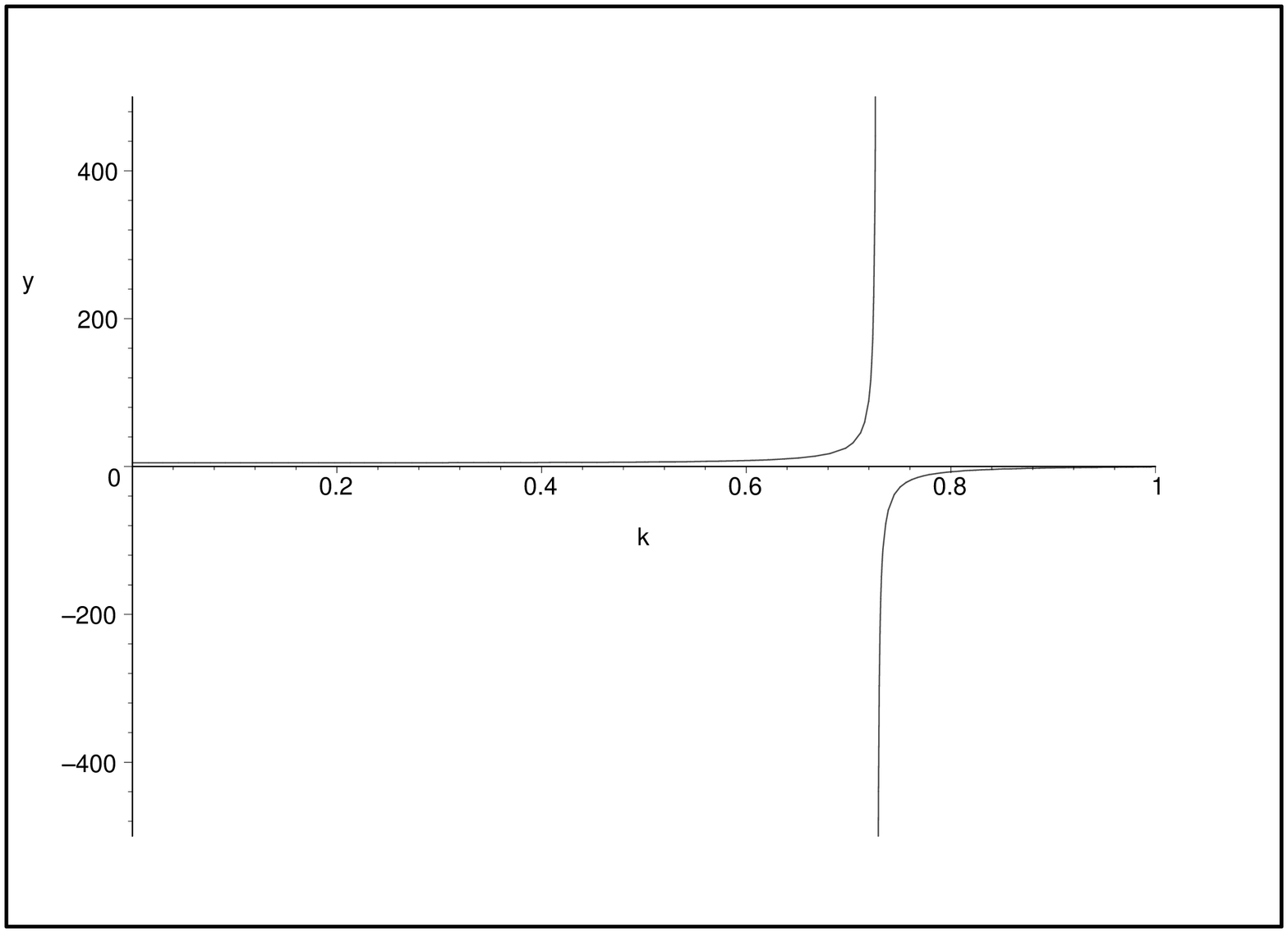}\quad 
\includegraphics[scale=0.29, angle=-90]{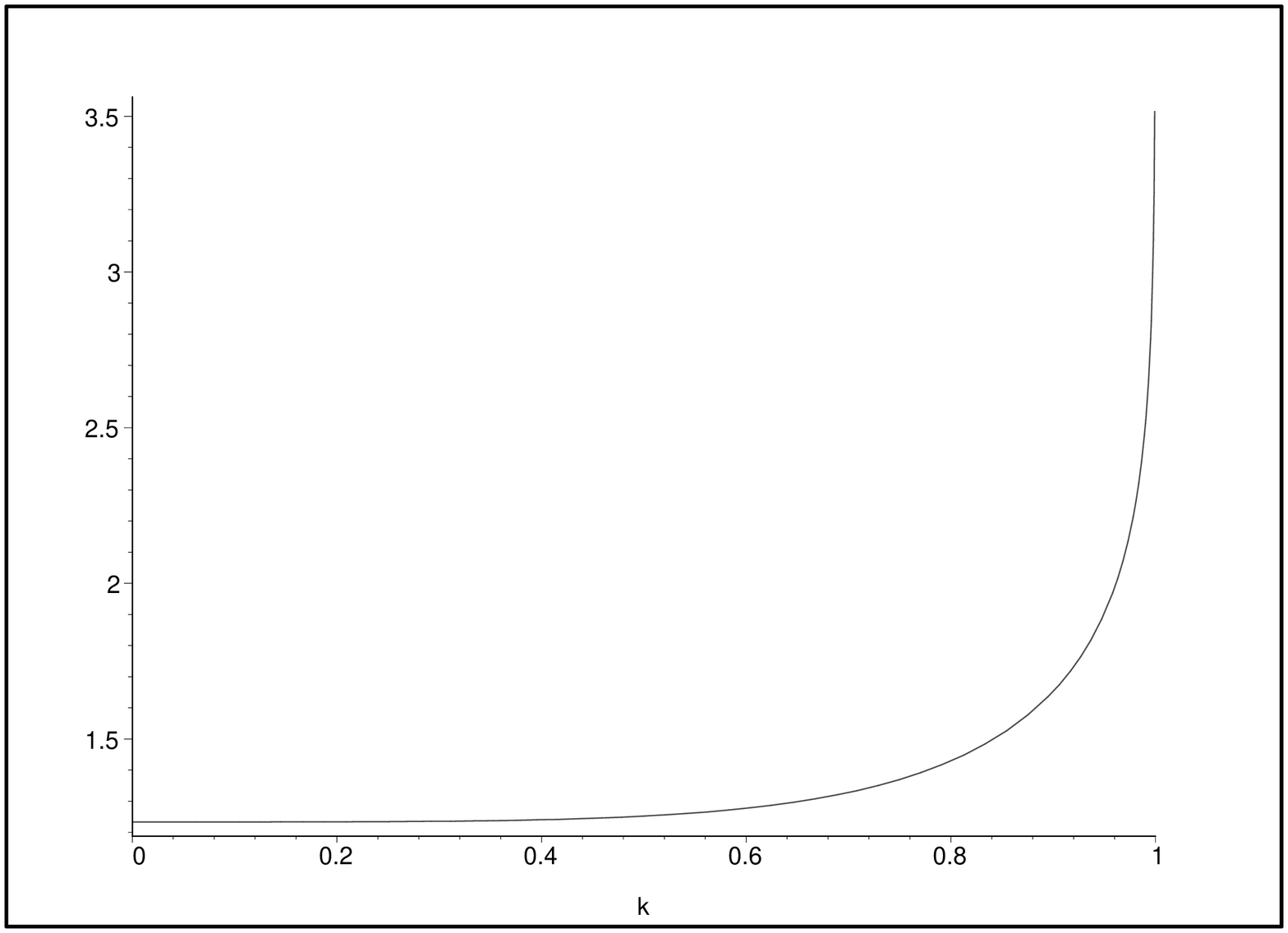}

\vspace{0.5cm}
\centerline{{\small Fig. 1: \textit{Graphic of $c(k)$  with } $L=7$ \hspace{1cm} Fig. 2: \textit{Graphic of $\tilde{a}(k)$ with}  $L=8$}}
\vspace{0.5cm}
Next, by making a similar analysis such as in the case of KdV equation in \cite{anguloNatali} ( see also \cite{angulo5}), we can obtain a smooth curve of positive cnoidal waves with the form $(\ref{curvaBBM})$, $c \in \left(c^{\ast}, +\infty\right)\longmapsto \phi_c \in H_{per}^n([0,L])$ for all $n\in\mathbb{N}$, such that $k:=k(c)$ is a strictly increasing smooth function of $c$. Moreover, we can determine that for $k \in (0,k_L)$ there is a unique $c \in \left(c^{\ast}, +\infty\right)$ such that $k(c)=k.$ Therefore, the function $w(k)$ defined above can be expressed as a function of $c$, $w=w(k(c))$, and it is a strictly increasing function (we will prove it later).\\

Now, since $\frac{K(k)}{K(k')} \in (0,L),$ for all $k\in(0,k_0)$, it follows that for  $c \in \left(c^{\ast}, +\infty\right)$ we obtain $w(k(c)) \in (1, +\infty)$. Therefore the mapping $c \in \left(c^{\ast}, +\infty\right)\longmapsto \psi_{w(k(c))} \in H_{per}^n([0,L])$ is a smooth curve for all $n \in \mathbb{N}$.\\

Concerning the well-posedness it is well-known that the BBM equation is local well-posed in $H^s_{per}$ for $s>1/2$ and global for $s\geq 1,$ which is enough for establish our theory of stability on the periodic case, but it is worth to note that Chen in \cite{hongqiu1} improved this result proving global well-posedness in $L^2_{per}.$\\

Next we present the stability result for the BBM equation.
\begin{theo}
Assume $L>2\pi$  fixed. If  $c>1+\frac{4\pi^2}{L^2-4\pi^2},$ then the periodic travelling wave solution $\phi_c$ in (\ref{curvaBBM}) is stable by the flow of the BBM equation.
\end{theo}
\proof
Note that $\phi_c = a{(k(c))}-\frac{24c}{L}\sqrt{\frac{w-1}{w}} +\frac{2c}{w}\psi_{w(k(c))}$, where 
\[a(k)=\frac{16cK}{L^2}\left[3E-(1+{k'}^2)K \right]+c-1.\]
Thus, for $s(k(c)):= a(k(c))- \frac{24c}{L}\sqrt{\frac{w-1}{w}}$, we can write $\phi_c(x)=s(k(c)) + \frac{2c}{w}\psi_{w(k(c))}(x)$. Therefore, we obtain easily that the Fourier coefficients of $\phi_c$ are
\[\widehat{\phi}_c(n)=\left\{
\begin{array}{l}
 a(k), \ \ \ n=0\\
\frac{12c \pi}{L^2}\ n \ \text{csch}\left(\sqrt{\frac{w}{w-1}}\frac{\pi n}{L}\right); \ \ \ n\neq 0,\  n \in \mathbb{Z}.
\end{array}\right.\]

It is worth to note that
\[a=\frac 1{L}\int_0^L\phi_c \ dx=s+\frac{2c}{wL}\int_0^L\psi_w\ dx,\]
then $\frac{1}{L}\int_0^L\psi_w\ dx=\frac{a-s}{2cw},$ that is, $\psi_w$ has average $\frac{a-s}{2cw}.$\\

Now, using the fact that $\frac{c-1}{c}=\frac{16K^2\sqrt{1-k^2+k^4}}{L^2}$ we can write $s$ as \[s(k)=c \left[\frac{16K^2}{L^2}\left( \sqrt{1-k^2+k^4}-2+k^2+3\frac{E}{K}\right)-\frac{24}{L^2}\frac{K(k)}{K(k')}\right]=: c\widetilde{s}(k).\]
The function $\widetilde{s}$ is a positive function in $(0,1)$, in particular in $(0,k_L),$ which does not have any root on the extremes of the interval $(0, k_L)$. We can also determine that $a(k)$ is a positive strictly increasing function (we wil prove it later). Therefore, we can follow the same ideas in  \cite{anguloNatali} to get that $\widehat{\phi}_c \in PF(2)$ discrete.\\

Next, it is easy to see that $\chi=-\frac{d}{dc}\phi_c$  satisfies  $\mathcal{L}\chi=\phi_c-\phi_c''$. Then by  Parseval theorem, it follows that $I=-\frac{L}{2}\frac{d}{dc}\left(\|(1+|\cdot|^2)^{\frac{1}{2}}\widehat{\phi}_c\|_{l^2}^2\right)$. But,
\begin{align*}
&\frac{d}{dc}\left(\|(1+|\cdot|^2)^{\frac{1}{2}}\widehat{\phi}_c\|_{l^2}^2\right)=2a(k)\frac{da}{dk}\frac{dk}{dc}
+C_1\sum_{\substack{n\in\mathbb{Z}\\n\neq 0}}(1+|n|^2)n^2 \text{csch}^2\left(\sqrt{\tfrac{w}{w-1}}\tfrac{\pi n}{L}\right)\\
&+C_2\left((w-1)^3w)\right)^{-1/2}\frac{dw}{dk}\frac{dk}{dc}\sum_{\substack{n\in\mathbb{Z}\\n\neq 0}}(1+|n|^2)n^3 \text{csch}^2\left(\sqrt{\tfrac{w}{w-1}}\tfrac{\pi n}{L}\right)\text{coth}\left(\sqrt{\tfrac{w}{w-1}}\tfrac{\pi n}{L}\right),
\end{align*}
where $C_1= C_1(L,c)>0$ and $C_2= C_2(L,c)>0$. Now, to prove that $I<0$ we only need to show that $\frac{dw}{dk}$ and $\frac{da}{dk}$ are positive, because $k=k(c)$ is a strictly increasing function and $b_n= (1+|n|^2)n^3 \ \text{csch}^2\left(\sqrt{\frac{w}{w-1}}\frac{\pi n}{L}\right)\text{coth}\left(\sqrt{\frac{w}{w-1}}\frac{\pi n}{L}\right)$ is clearly a positive sequence. Hence, we have.
\[\frac{dw}{dk}=\frac{2L^2K'K\left[k'\frac{dK}{dk}-K\frac{dK'}{dk}\right]}{(L^2{K'}^2-K^2)^2}, \ \ \ \ \forall k\in(0,k_0).\]
Since $\frac{dK}{dk}>0$ and $\frac{dK'}{dk} <0$ we get that $\frac{dw}{dk}>0$ for all $k\in(0,k_0)$.\\

Finally, note that
\begin{align*}
a(k)&=c\left\{\frac{16K}{L^2}\left[3E -(2-k^2)K\right]+ \frac{c-1}{c}\right\}\\
&=\frac{16K^2c}{L^2}\left[3\frac{E}{K}-2+k^2+\sqrt{1-k^2+k^4}\right]=: c\widetilde{a}(k).
\end{align*}
Since $\widetilde{a}(k)$ and $c(k)$ are an increasing positive function (see Figure 1 and 2) we have that $\frac{da}{dk} > 0$. Therefore, $I<0$ and the positive cnoidal waves $\phi_c$ are stable in $H_{per}^1([0,L])$ by the periodic flow of the BBM equation.
\endproof

\begin{remark}
In an upcoming work we will present the stability theory for the modified BBM equation $(p=2)$ and the critical BBM equation ($p=4).$ 
\end{remark}

\textbf{Acknowledgements:} The first author was supported by CAPES/Brazil Grant. The third author was  supported by CNPq/Brazil doctoral fellowship at State University  of Campinas.


\end{document}